\documentclass[11pt,a4paper,reqno]{amsart}
\usepackage{amsthm,amsmath,amsfonts,amssymb,amsxtra,appendix,bookmark,dsfont,latexsym,epsfig,graphicx,stix}

\makeatletter
\usepackage{hyperref}
\usepackage{delarray,color}
\usepackage{euscript}

\usepackage{hyperref}

\makeatletter
\def\@tocline#1#2#3#4#5#6#7{\relax
  \ifnum #1>\c@tocdepth 
  \else
    \par \addpenalty\@secpenalty\addvspace{#2}%
    \begingroup \hyphenpenalty\@M
    \@ifempty{#4}{%
      \@tempdima\csname r@tocindent\number#1\endcsname\relax
    }{%
      \@tempdima#4\relax
    }%
    \parindent\z@ \leftskip#3\relax \advance\leftskip\@tempdima\relax
    \rightskip\@pnumwidth plus4em \parfillskip-\@pnumwidth
    #5\leavevmode\hskip-\@tempdima
      \ifcase #1
       \or\or \hskip 1em \or \hskip 2em \else \hskip 3em \fi%
      #6\nobreak\relax
      \dotfill
      \hbox to\@pnumwidth{\@tocpagenum{#7}}
    \par
    \nobreak
    \endgroup
  \fi}
\makeatother

\newcommand{\bdm}{\begin{displaymath}}
\newcommand{\edm}{\end{displaymath}}
\newcommand{\bdn}{\begin{eqnarray}}
\newcommand{\edn}{\end{eqnarray}}
\newcommand{\bay}{\begin{array}{c}}
\newcommand{\eay}{\end{array}}
\newcommand{\ben}{\begin{enumerate}}
\newcommand{\een}{\end{enumerate}}

\newcommand{\N}{\mathbb{N}}
\newcommand{\R}{\mathbb{R}}
\newcommand{\C}{\mathbb{C}}
\newcommand{\cW}{\mathcal{W}}

\newcommand{\PsiLau}{\Psi_{\rm Lau}}
\newcommand{\PsiLaun}{\Psi_{\rm Lau} ^{(\ell)}}
\newcommand{\cLau}{c _{\rm Lau}}

\newcommand{\cH}{\mathcal{H}}

\newcommand{\TFM}{\mathcal{M} ^{\rm TF}}
\newcommand{\ETF}{\mathcal{E} ^{\rm TF}}

\newcommand{\TFmin}{\sigma ^{\rm TF}}
\newcommand{\TFpot}{\Phi ^{\rm TF}}

\def\j{\mathbf{E}}
\def\bam{\overline{\mathcal{A}}_\rho}

\newtheorem{defi}{Definition}[section]

\newtheorem{theorem}{Theorem}[section]

\newtheorem{definition}[theorem]{Definition}

\theoremstyle{remark}
\newtheorem{remark}[theorem]{Remark}

\newcommand{\beq}{\begin{equation}}
\newcommand{\eeq}{\end{equation}}
\newcommand{\one}{{\ensuremath {\mathds 1} }}

\newcommand{\rhoMF}{\rho^{\rm MF}}
\newcommand{\EMF}{F^{\rm MF}}
\newcommand{\cEMF}{\mathcal{F}^{\rm MF}}

\numberwithin{equation}{section}

\newcommand{\bx}{\mathbf{x}}
\newcommand{\by}{\mathbf{y}}
\newcommand{\im}{\mathrm{i}}
\newcommand{\1}{\mathds{1}}
\newcommand{\mub}{\boldsymbol{\mu}}
\newcommand{\nub}{\boldsymbol{\nu}}
\newcommand{\cZ}{\mathcal{Z}}
\newcommand{\cF}{\mathcal{F}}

\newcommand{\coul}{\mathrm{\mathbf{C}}}

\def\Xint#1{\mathchoice
   {\XXint\displaystyle\textstyle{#1}}%
   {\XXint\textstyle\scriptstyle{#1}}%
   {\XXint\scriptstyle\scriptscriptstyle{#1}}%
   {\XXint\scriptscriptstyle\scriptscriptstyle{#1}}%
   \!\int}
\def\XXint#1#2#3{{\setbox0=\hbox{$#1{#2#3}{\int}$}
     \vcenter{\hbox{$#2#3$}}\kern-.5\wd0}}
\def\dashint{\Xint-}

\def\W{\mathcal{W}}

\begin{document}

\title{The classical Jellium and the Laughlin phase}

\author[N. Rougerie]{Nicolas Rougerie}
\address{Ecole Normale Sup\'erieure de Lyon \& CNRS,  UMPA (UMR 5669)}
\email{nicolas.rougerie@ens-lyon.fr}





\date{March, 2022}

\begin{abstract}
I discuss results bearing on a variational problem of a new type, inspired by fractional quantum Hall physics. In the latter context, the main result reviewed herein can be spelled as ``the phase of independent quasi-holes generated from Laughlin's wave-function is stable against external potentials and weak long-range interactions''. The main ingredient of the proof is a connection between fractional quantum Hall wave-functions and statistical mechanics problems that generalize the 2D one-component plasma (jellium model). Universal bounds on the density of such systems, coined ``Incompressibility estimates'' are obtained via the construction of screening regions for any configuration of points with positive electric charges. The latter regions are  patches of constant, negative electric charge density, whose shape is optimized for the total system (points plus patch) not to generate any electric potential in its exterior. 
\end{abstract}

\maketitle

\begin{center}
 \emph{Dedicated to Elliott H. Lieb, on the occasion of his 90th birthday. With admiration.}
\end{center}

\tableofcontents

\section{Foreword}

The purpose of this note is to discuss a somewhat unusual variational problem in quantum mechanics (introduced and investigated in~\cite{RouSerYng-13a,RouSerYng-13b,RouYng-14,RouYng-15,RouYng-17,RouYng-19,LieRouYng-16,LieRouYng-17,OlgRou-19}) and its relation to the statistical mechanics of Coulomb systems, in particular to recent results bearing on the latter. The mathematical set-up is first described as concisely as possible in Section~\ref{sec:prob}. Next, the physical motivation for such investigations and the elements of context allowing to interpret the main result are discussed in~Section~\ref{sec:FQHE}. The keywords here are the fractional quantum Hall effect and the Laughlin function.

As regards proofs, everything proceeds from the ``plasma analogy'', which maps trial states for the original, quantum, variational problem onto classical statistical mechanics Hamiltonians. More precisely, the modulus squared of the quantum wave-functions have an interpretation in terms of Gibbs states of effective 2D Coulomb Hamilton functions. In Section~\ref{sec:jellium} we give a brief review on the statistical mechanics of the simplest of such equilibria: the  jellium (or one-component plasma), namely a system of classical point charges interacting via repulsive Coulomb forces and attracted to a neutralizing background of opposite charge. Basic questions (for which the recent~\cite{Lewin-22} provides an in-depth review) concern 
\begin{itemize}
 \item the existence of the thermodynamic limit for homogeneous systems, as investigated first in pioneering works by Elliott H. Lieb and co-workers~\cite{LieNar-76,LebLie-69,LieLeb-72} and then taken up e.g. in~\cite{SarMer-76,Fefferman-85,HaiLewSol_proc-08}.
\item the local density approximation for inhomogeneous systems, investigated in a long series of works by Serfaty and co-workers, reviewed e.g. in~\cite{Serfaty-14,Serfaty-15,Serfaty-17}.
\end{itemize}
In our applications to the variational problem to be described shortly, our needs are somewhat more specific. What we need are ``incompressibility estimates'': universal local density upper bounds for a certain class of generalized Coulomb systems. The simplest of such estimates follows from an unpublished theorem of Elliott: in the ground state of a classical jellium, the minimal distance between two point charges is bounded below, uniformly in the thermodynamic limit. The reason is that each point charge neutralizes a circular patch of the background around it, creating a region where it is never energetically favorable to put another point charge.

This method was generously communicated by Elliott to Sylvia Serfaty in private conversation. Generalizations thereof played a key role in the aforementioned program~\cite{RouSer-14,RotSer-14,PetSer-14}. Passed on to the present author and Jakob Yngvason, it allowed to prove the first unconditional incompressibility estimate~\cite{RouYng-15}, after the problem had been isolated and first progress made in~\cite{RouYng-14}. To obtain a satisfying bound, whose use unlocked the main theorem from Section~\ref{sec:prob} below, a new key idea was necessary, proposed by Elliott to Jakob Yngvason and the present author. One can in fact associate to \emph{any} set of point charges a patch of the background such that the potential generated by the ensemble vanishes outside of the patch\footnote{During the writing of this text I became aware of the fact that this concept predates our work~\cite{GusSha-95,GusPut-07,Sakai-82,Gustafsson-02}. See Remark~\ref{rem:balayage} below.}. It is never energetically favorable to add another charge inside such a patch. The construction of such \emph{screening regions} is not straightforward, nor is their use to complete the proof of the needed incompressibility estimates. All this is the topic of~\cite{LieRouYng-16,LieRouYng-17}. Further developments are in~\cite{RouYng-17,OlgRou-19}, as we will explain below (see~\cite{Rougerie-xedp19} for another exposition).

Before proceeding, some directions related to this note, but whose discussion goes beyond its scope, are worth mentioning. Indeed, Elliott made other inspiring contributions to the study of the Laughlin function~\cite{JanLieSei-08,JanLieSei-09} or the statistical mechanics of jellium-like models~\cite{LewLieSei-18,LewLieSei-19a,LewLieSei-19b}. 

\section{A unusual variational problem}\label{sec:prob}

\subsection{Statements}

We start from a very standard Hamilton function 
\begin{equation}\label{eq:start Hamil}
 \R^{2N} \ni \left( \bx_1,\ldots,\bx_N \right) \mapsto \sum_{j=1}^N V(\bx_j) + \lambda \sum_{1\leq i < j \leq N} W (\bx_i - \bx_j) 
\end{equation}
where $V,W:\R^2 \mapsto \R$ are respectively a one-body and a two-body potential and $\lambda\in \R$ is a coupling constant. 

We are interested in minimizing the expectation value of the above in particular quantum wave-functions of the following form. Let $B>0$ and 
\begin{equation}\label{eq:PsiLau}
 \PsiLau ^{(\ell)} (z_1,\ldots,z_N):= \cLau \prod_{1\leq i < j \leq N} (z_i-z_j) ^\ell e^{-\frac{B}{4}\sum_{j=1} ^N |z_j| ^2}
\end{equation}
be the Laughlin function of exponent $\ell \in \N^*$, where the planar coordinates $\bx_1,\ldots,\bx_N $ are identified with complex numbers $z_1,\ldots,z_N$ and $\cLau=\cLau (\ell)$ is a $L^2$-normalization constant. For any $F:\C^N \mapsto \C$ analytic in all its argument and symmetric in the sense that 
\begin{equation}\label{eq:sym}
F (z_{\sigma(1)},\ldots,z_{\sigma(N)}) =  F (z_1,\ldots,z_N) 
\end{equation}
for any permutation $\sigma$ of $N$ indices, let  
\begin{equation}\label{eq:PsiF}
 \Psi_F (z_1,\ldots,z_N):= c_F \PsiLaun (z_1,\ldots,z_N) F (z_1,\ldots,z_N)
\end{equation}
with $c_F$ a $L^2$-normalization constant.   

We shall discuss the following problem 
\begin{equation}  \label{eq:qm_energy}
E (N,\lambda)=\inf\Big\{\mathcal{E}_{N,\lambda}[\Psi_F]\;|\;\Psi_F \mbox{ of the form~\eqref{eq:PsiF}},\,\int_{\mathbb{R}^{2N}}|\Psi_F|^2=1\Big\}
\end{equation}
where 
\begin{equation} \label{eq:many_body_energy}
\mathcal{E}_{N,\lambda}[\Psi_F]=\Big\langle\Psi_F\Big|\sum_{j=1}^N V(x_j)+\lambda \sum_{i<j}W(x_i-x_j)\Big|\Psi_F\Big\rangle_{L^2}.
\end{equation}
What we aim at is a significant simplification in the $N\to \infty$ limit. Namely, consider the simplest functions of the form~\eqref{eq:PsiF}:
\begin{equation}\label{eq:Psif}
 \Psi_f (z_1,\ldots,z_N):= c_f \PsiLaun(z_1,\ldots,z_N) \prod_{j=1} ^N f(z_j)
\end{equation}
where $f:\C \mapsto \C$ is analytic and $c_f$ is a normalization constant. Define a restricted infimum by setting 
\begin{equation} \label{eq:qh_energy}
e (N,\lambda)=\inf\Big\{\mathcal{E}_{N,\lambda}[\Psi_f]\;|\;\Psi_f \text{ of the form \eqref{eq:Psif}},\,\int_{\mathbb{R}^{2N}}|\Psi_f|^2=1\Big\}.
\end{equation}
Obviously $E(N,\lambda) \leq e(N,\lambda)$. What we would like to prove is that 
\begin{equation} \label{eq:anticipation}
\boxed{E (N,\lambda)\simeq e (N,\lambda)\quad\text{as}\;N\to\infty \mbox{ with } \lambda \mbox{ fixed}.}
\end{equation}
In fact, functions from our variational set~\eqref{eq:PsiF} naturally live over thermodynamically large length scales $\sim \sqrt{N}$.  It is hence natural to scale the potentials $V$ and $W$ accordingly. We thus set, for fixed functions $v,w$,
\begin{equation}\label{eq:rescaled_V}
V (x) = v\left(N^{-1/2} x\right)
\end{equation}
and (the $N^{-1}$ pre-factor ensures that the potential and interaction energies stay of the same order when $N\to \infty$)
\begin{equation} \label{eq:rescaled_w}
W (x) = N ^{-1} w \left(N^{-1/2} x\right).
\end{equation}
We are now ready to state the main result we want to discuss. It was proved in~\cite{OlgRou-19} after the (simpler but still highly non-trivial) $\lambda = 0$ version was obtained in~\cite{LieRouYng-17,RouYng-17}. We do not state precise or optimal assumptions, nor associated corollaries, for the sake of a simpler exposition.

\begin{theorem}[\textbf{Energy of the Laughlin phase}]\label{thm:ener}\mbox{}\\
Assume that $v$ and $w$ are smooth fixed functions. Assume that $v$ goes to $+\infty$ polynomially at infinity, and that it has finitely many non-degenerate critical points. There exists $\lambda_0 >0$ such that
$$ \frac{E(N,\lambda)}{e(N,\lambda)} \underset{N\to \infty}{\to} 1$$
with $B > 0$ fixed, $\ell>0$ a fixed integer and $|\lambda| \leq \lambda_0$. 
\end{theorem}

\subsection{Proof outline and the connection with jellium}\label{sec:analogy}

To set the stage for the material to be discussed below, we briefly sketch the main steps of the proof of Theorem~\ref{thm:ener}. 

\medskip

\noindent\textbf{Plasma analogy}. The main difficulty is to understand what the densities of wave-functions of the form~\eqref{eq:PsiF} have in common. We use Laughlin's plasma analogy from~\cite{Laughlin-83,Laughlin-87}, writing $|\Psi_F|^2$ as a Boltzmann-Gibbs factor,
\begin{equation}\label{eq:plasma}
 |\Psi_F (z_1,\ldots,z_N)| ^2 = \frac{1}{\mathcal{Z}_F} \exp\left( - H_F (z_1,\ldots,z_N)\right)
\end{equation}
with $\mathcal{Z}_F$ ensuring $L^1$-normalization (partition function of the effective plasma) and $H_F$ an effective Hamilton function
\begin{equation}\label{eq:class hamil}
H_F (z_1,\ldots,z_N) = \frac{B}{2} \sum_{j=1} ^N |z_j| ^2 - 2 \ell \sum_{1\leq i < j \leq N} \log |z_i-z_j| - 2 \log \left| F (z_1,\ldots,z_N) \right|. 
\end{equation}
Hence $|\Psi_F|^2$ is the probability density of particles minimizing the classical free energy associated with $H_F$, that is, realizing the infimum 
$$ -\log \mathcal{Z}_F = \inf_{\boldsymbol{\mu}_N} \left\{ \int_{\R^{2N}} H_F \boldsymbol{\mu}_N + \int_{\R^{2N}} \boldsymbol{\mu}_N \log \boldsymbol{\mu}_N\right\}$$
over probability measures on $\R^{2N}$. 

This rewriting is fruitful because the latter functional has an interpretation in terms of 2D electrostatics. The potential $\Phi$ generated by a charge distribution $\sigma$ is given by 
$$ - \Delta \Phi = \sigma, \quad \Phi = -\frac{1}{2\pi}\log|\,.\,| \star \sigma.$$
Hence $H_F$ is the energy of $N$ mobile negatively charged 2D particles (at locations $z_1,\ldots,z_N \in \C \leftrightarrow \R^2$) of charge $- \sqrt{4\pi \ell}$

\medskip 

\noindent \textbf{1}. Interacting among themselves via repulsive Coulomb forces. 

\medskip
 
\noindent \textbf{2}. Attracted to a fixed uniform background of positive charge density 
$$\Delta \left(\frac{B}{4\sqrt{\pi \ell}} |z|^2\right) = \frac{B}{\sqrt{\pi \ell}}.$$

\medskip

\noindent \textbf{3}. Feeling the potential 
\begin{equation}\label{eq:crazy pot}
\cW := -2 \log |F| 
\end{equation}
generated by additional ``phantom'' positive charges. The location of the latter can be essentially arbitrary, and correlated with the positions of $z_1,\ldots,z_N$, but their charge must be positive because 
\begin{equation}\label{eq:superharm}
 - \Delta_{z_j} \cW  \geq 0 
\end{equation}
for any $j$ (recall that $F$ is analytic).

\medskip

The last equation~\eqref{eq:superharm} is key to the method I expose here. For the more specific functions~\eqref{eq:Psif} the interpretation can be made cleaner. Taking $f$ to be a polynomial 
\begin{equation}\label{eq:polynomial}
 f(z) = \prod_{k=1}^K (z-a_k) 
\end{equation}
with zeros $a_1,\ldots,a_K \in \C$ we have that 
\begin{equation}\label{eq:qhonebod}
 \cW = -2 \sum_{k=1}^K \sum_{j=1}^N \log |z_j-a_k| 
\end{equation}
corresponds to the electrostatic interaction of our mobiles charges $z_1,\ldots,z_N$ with fixed point charges $a_1,\ldots,a_K$ of the opposite sign.

\medskip

\noindent\textbf{Flocking energy and the main strategy.} There is a third energy that we use as an intermediary step to relate $E(N,\lambda)$ to $e(N,\lambda)$. This is the flocking energy 
\begin{equation}\label{eq:flocking}
E^{\rm flo} (N,\lambda) := \inf \left\{\int_{\R^2} V \varrho + \frac{\lambda}{2} \iint_{\R^2} \varrho(\bx) W (\bx-\by) \varrho(y) d\bx d\by  , \: 0 \leq \varrho \leq \frac{B}{2\pi \ell}, \: \int_{\R^2} \varrho = N \right\}. 
\end{equation}
Note the mean-field character: we minimize over a single one-body density, assuming negligible correlations. However, the original nature of the problem is retained in the upper constraint  
\begin{equation}\label{eq:flo bound}
 \varrho \leq \frac{B}{2\pi \ell} 
\end{equation}
imposed on admissible densities. This is a kind of ``super-Pauli principle'': particles are prevented to gather in space beyond a certain fixed density, lower than the usual Pauli principle ($\ell=1$) would require. This is reminiscent of simple models~\cite{BurChoTop-15} for flocks or swarms of animals, e.g. birds, a direction where Elliott also made contributions~\cite{FraLie-16}. 

Since we clearly have $E(N,\lambda) \leq e (N,\lambda)$, the main theorem will follow from the following arguments:

\medskip

\noindent\underline{(1) Energy upper bound}. By constructing suitable trial states we have 
$$ e (N,\lambda) \lessapprox E^{\rm flo} (N,\lambda).$$
The idea is developed in~\cite{RouYng-17,OlgRou-19} after particular cases were dealt with in~\cite{RouSerYng-13a,RouSerYng-13b}. One writes the function $f$ in~\eqref{eq:Psif} as in~\eqref{eq:polynomial} and optimizes the number and locations of the zeroes $a_1,\ldots,a_K$ for the effective plasma to have a density matching that of the solution to~\eqref{eq:flocking}. For this to be possible it is crucial that the latter saturates the density upper bound~\eqref{eq:flo bound} everywhere on its support, i.e. is in the ``solid'' phase~\cite{FraLie-16}. This is where we use the assumption that $\lambda$ is small enough. For reasons explained in more details in~\cite{OlgRou-19}, such a constraint is a natural requirement. One should \emph{not} expect Theorem~\ref{thm:ener} to hold for repulsive interaction potentials if $\lambda$ is allowed to be arbitrarily large. 

\medskip

\noindent\underline{(2) Energy lower bound, $\lambda=0$.} With no interaction potential in~\eqref{eq:start Hamil}, the flocking energy is a simple bath-tub problem~\cite[Theorem~1.14]{LieLos-01} and the many-body problem~\eqref{eq:many_body_energy} only depends on the one-particle density 
\begin{equation}\label{eq:density}
\varrho_{\Psi_F} (\bx) := N \int_{\R^{2(N-1)}} \left| \Psi_F (\bx,\bx_2,\ldots,\bx_N) \right|^2 d\bx_2 \ldots d\bx_N. 
\end{equation}
What is needed is a proof that the set of admissible $\varrho_F$ is in some sense included in the variational set defining~\eqref{eq:flocking}, namely 
\begin{equation}\label{eq:incomp vague}
\varrho_F \lessapprox \frac{B}{2\pi \ell}  
\end{equation}
in some appropriate sense, for \emph{any} symmetric analytic function $F$ entering~\eqref{eq:PsiF}. This is what we called an ``incompressibility estimate'' in~\cite{RouYng-14,RouYng-15}. The lower bound 
$$ E (N,0) \gtrapprox E^{\rm flo} (N,0)$$
essentially follows.

\medskip

\noindent\underline{(3) Energy lower bound, $\lambda\neq 0$.} With interactions there is an extra mean-field limit to deal with (recall the scaling~\eqref{eq:rescaled_V}-\eqref{eq:rescaled_w}), and it is necessary to understand why correlations can be neglected. In~\cite{OlgRou-19} we used the approach to classical mean-field limits based on the de Finetti-Hewitt-Savage theorem, see~\cite[Chapter~2]{Rougerie-LMU,Rougerie-cdf} for review. The energy now genuinely depends on the pair density
\begin{equation}\label{eq:pair density}
\varrho_{\Psi_F}^{(2)} (\bx,\by) := N (N-1) \int_{\R^{2(N-2)}} \left| \Psi_F (\bx,\by,\bx_3,\ldots,\bx_N) \right|^2 d\bx_3 \ldots d\bx_N 
\end{equation}
and we prove that we can approximate, for any $F$,
\begin{equation}\label{eq:deF}
\varrho_{\Psi_F}^{(2)} (\bx,\by) \simeq \int_{\varrho \leq \frac{B}{2\pi\ell}} \varrho (\bx) \varrho(\by) dP(\varrho)
\end{equation}
where $P$ is a probability measure over one-particle densities satisfying 
$$ 0 \leq \varrho \leq \frac{B}{2\pi\ell}, \int_{\R^2} \varrho = N.$$
The most important is that $P$ only charges functions satisfying the incompressibility bound~\eqref{eq:incomp vague}. 

\medskip

To obtain an energy upper bound at $\lambda \neq 0$ we also prove that, for the pair density associated to a trial state $\Psi_f$ of the form~\eqref{eq:Psif}
$$ \varrho_{\Psi_f} ^{(2)}(\bx,\by) \simeq \varrho_{\Psi_f} (\bx) \varrho_{\Psi_f} (\by),$$
which is again a mean-field kind of problem, but for the effective plasma Hamiltonians, not the original physical one. 

In this note we wish to focus on the main ingredients, namely the incompressibility estimates necessary for steps (2) and (3) above. After partial progress in~\cite{RouYng-14,RouYng-15},  together with Elliott and Jakob Yngvason we obtained~\cite{LieRouYng-16,LieRouYng-17} the following version of~\eqref{eq:incomp vague}:

\begin{theorem}[\textbf{Incompressibility estimates in the mean}]\label{thm:incomp}\mbox{}\\ 
For any\footnote{The result is obtained in~\cite{LieRouYng-17} for any $\alpha >1/4$, but the proof requires mild additional trapping conditions on $\Psi_F$, see~\cite[Remark~4.2]{OlgRou-19}.} $\alpha>(\sqrt{5}-1)/4$, any disk  $D$ of radius $N^{\alpha}$ and any (sequence of) states $\Psi_F$ of the form~\eqref{eq:PsiF} we have 
\begin{equation}\label{bound}
\int_{D} \varrho_{\Psi_F} \leq \frac{B}{2\pi\ell} |D| (1+o(1))
\end{equation}
where $|D|$ is the area of the disk and $o(1)$ tends to zero as $N\to \infty$. 
\end{theorem}

The above means that (with additional mild assumptions,~\cite[Section~5.1]{LieRouYng-17}) the desired bound holds in the sense of averages (on disks, or actually on any other nice set) of length-scale $L\gg N^{1/4}$. Note that the thermodynamic length scale is $L\sim N^{1/2}$ in this problem, and the typical inter-particle distance $L\sim B^{-1/2}$ (the magnetic length, fixed in our convention). We expect that the bound actually holds on any mesoscopic length scale $L\gg 1$. 

To obtain the stronger notion of incompressibility~\eqref{eq:deF}, we construct the measure $P$ in the standard de Finetti-Hewitt-Savage-Diaconis-Freedman way~\cite{Rougerie-cdf,Rougerie-LMU}. The difficulty is to prove that it charges only densities satisfying the appropriate bound, which follows from the following result from~\cite{OlgRou-19}:

\begin{theorem}[\textbf{Probability of violating the incompressibility estimate}] \label{thm:exponential_bound}\mbox{}\\
Let $D$ be a disk of radius $N^\alpha$ as above. Let $|\Psi_F|^2$ be the probability measure on $\R^{2N}$ associated with $\Psi_F$ of the form~\eqref{eq:PsiF}. Denote $\mathbb{P}_F (A)$ the associated probability of events $A\subset \R^{2N}$ and $\sharp(A)$ the cardinal of a discrete set $A$. Then, for any $\varepsilon>0$,
\begin{equation} \label{eq:exponential_bound}
\mathbb{P}_F \left(\left\{X_N \in \mathbb{R}^{2N} \mbox{ such that } \sharp(X_N \cap D) > (1+\varepsilon)\frac{B|D|}{2\pi\ell} \right\}\right)\le \exp\left(-C\varepsilon N^{\sqrt{5}-1}\right).
\end{equation}
\end{theorem}

\noindent This ``large deviation bound'' implies that not only~\eqref{eq:incomp vague} but also the stronger bound
\begin{equation}\label{eq:incomp vague plus}
\varrho_{\Psi_F} ^{(k)} \lessapprox \left(\frac{B}{2\pi \ell}\right)^k  
\end{equation}
holds for all $k$-particles reduced densities
\begin{equation}\label{eq:red density}
\varrho_{\Psi_F}^{(k)} (\bx_1,\ldots,\bx_k) := \frac{N!}{(N-k)!} \int_{\R^{2(N-k)}} \left| \Psi_F (\bx_1,\ldots,\bx_k,\bx_{k+1},\ldots,\bx_N) \right|^2 d\bx_{k+1} \ldots d\bx_N 
\end{equation}
at least if $k$ is fixed in the limit $N\to \infty$. The de Finetti-Hewitt-Savage-Diaconis-Freedman theorem and arguments from~\cite{FouLewSol-15} (originally used for fermionic semi-classical measures on phase-space) then imply~\eqref{eq:deF}. 

\section{Physical context: the fractional quantum Hall effect}\label{sec:FQHE}

We now explain how the above strange variational problem sheds light on some aspects of fractional quantum Hall physics. This section can be skipped by readers who prefer to jump to the connection between the incompressibility bounds of Theorems~\ref{thm:incomp}-\ref{thm:exponential_bound} and the statistical mechanics of the one-component plasma.

\bigskip

\noindent\textbf{Summary}. It is useful to first set a road-map for this section: 
\begin{itemize}
 \item In a typical fractional quantum Hall experiment at filling factor $\nu \sim \ell ^{-1}$, the Laughlin function~\eqref{eq:PsiLau} is a well-educated guess for the system's vacuum. It freezes the magnetic kinetic energy and reduces a lot the short-range part of the interaction. 
 \item The quasi-holes wave-functions~\eqref{eq:Psif}-\eqref{eq:polynomial} are proposed to describe the vacuum's excitations in response to a slightly smaller filling factor, to impurities in the sample, to residual interactions, to external fields etc ...
 \item The arguments leading to the above points in fact allow in generality to work with the class of states~\eqref{eq:PsiF}. The further reduction to~\eqref{eq:Psif} is motivated by (legitimate !) arguments based on simplicity/guesses, and, ultimately, experimental confirmation.
 \item That the Laughlin phase of uncorrelated Laughlin quasi-holes does in fact emerge as the effective ground state of the system is justified in the thermodynamic limit by Theorem~\ref{thm:ener} and corollaries.   
\end{itemize}

The rest of the section is meant as a clarification of the above summary.

\bigskip

\noindent\textbf{The many-body quantum Hamiltonian}. We start from a basic Hamiltonian for the quantum 2D electron gas (in adimentionalised form $\hbar = c = e = 2m = 1$)
\begin{equation}\label{eq:mag hamil}
H_N^{\rm QM} = \sum_{j=1} ^N \left[\left( -\im \nabla_{\bx_j} - \frac{B}{2} \bx_j ^\perp \right) ^2 + V (\bx_j) \right]+ \sum_{1\leq i < j \leq N} W(\bx_i-\bx_j) 
\end{equation}
acting on $L_{\rm asym} ^2 (\R ^{2N}) $, the Hilbert space for $N$ 2D fermionic particles. Here $\bx^{\perp}$ denotes the vector $\bx\in \R^2$ rotated by $\pi/2$ counter-clockwise, so that 
$$ \mathrm{curl} \frac{B}{2} \bx ^\perp = B$$
and thus $\frac{B}{2} \bx ^\perp$ is the vector potential of a uniform magnetic field, expressed in symmetric gauge. In view of our choice of units, $B$ is actually $\sqrt{\alpha}$ times the physical magnetic field, with $\alpha = e^2 /(\hbar c) \sim 1/137$ the fine structure constant, see e.g.~\cite[Section~2.17]{LieSei-09}. 

We take into account an external potential $V:\R^2 \mapsto \R$ modeling trapping and/or impurities in the sample, and repulsive pair interactions $W:\R^2 \mapsto \R$ between particles. Typically $W$ should be the 3D Coulomb kernel (with $\alpha$ the fine structure constant again)
$$ W (\bx-\by) = \frac{\alpha}{|\bx-\by|}$$
or some screened version. We have made the customary assumption that the magnetic field is strong enough to polarize all the electrons' spins.

\bigskip

\noindent\textbf{The quantum Hall effect}~\cite{Jain-07,Girvin-04,Goerbig-09,StoTsuGos-99,Laughlin-99} is a peculiar feature of the transport properties of 2D electron gases under strong perpendicular magnetic fields. The main experimental findings (see Figure~$1$) are plateaux in the Hall (transverse) resistance $R_{xy}$ at particular quantized values, accompanied with huge drops in the longitudinal resistance $R_{xx}$. The extremely precise quantization\footnote{The general ``interpolating'' linear shape of the curve $B\mapsto R_{xy}$ is the \emph{classical} Hall effect, a neat way of measuring charge carrier densities in samples.} to particular values of $R_{xy}$ (read on the vertical axis of Figure~$1$) has an interpretation in terms of topological invariants of the system~\cite{BelSchEls-94,Frohlich-92,Frohlich-95}, but that is not what we focus on here. Instead, looking at the horizontal axis of Figure~$1$, we see that the particular features occur around special values (the numbers associated with arrows on the picture) of the filling factor of the system 
\begin{equation}\label{eq:filling}
\nu := \frac{hc}{e} \frac{\rho}{B} 
\end{equation}
with $\rho$ the electrons' density, $B$ the applied magnetic field and $h,c,e$ respectively Planck's constant, the speed of light and the elementary charge. In this note we (partially) address only the question ``\emph{why} does something special happen at these parameter values~?'' without touching much on the ``\emph{how} does the particular observed experimental signature emerge~?''

 \begin{figure}\label{fig:FQHE}
\includegraphics[width=10cm]{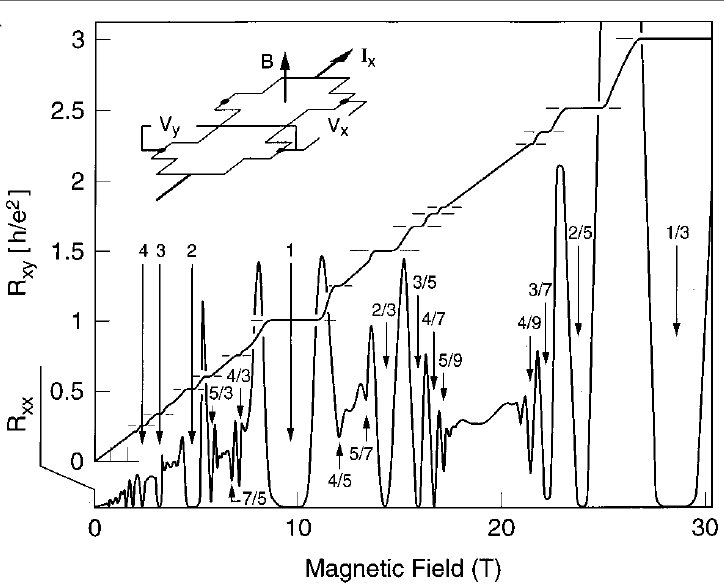}
\caption{The fractional quantum Hall effect~\cite{StoTsuGos-99}. Sketch of the experimental sample in top-left corner. Plots of the longitudinal $R_{xx} = V_x /I_x$ and transverse (Hall) $R_{xy} = V_y/I_x$ resistances as a function of the magnetic field.}
\end{figure}

\bigskip

\noindent\textbf{Landau levels.} The workhorse of the quantum Hall effect is the quantization of kinetic energy levels in the presence of a magnetic field. Namely, the appropriate kinetic energy operator for a 2D particle in a perpendicular magnetic field $B$ is 
\begin{equation}\label{eq:Landau}
 H = \left( -\im \nabla_{\bx} - \frac{B}{2} \bx ^\perp \right) ^2 
\end{equation}
acting on $L^2 (\R^d)$. 

The energy levels (eigenvalues) of the above are well-known~\cite{RouYng-19,Jain-07} to be $2B (n+1/2)$ for integer $n$, since one can write 
$$ H = 2B \left( a^\dagger a + \frac{1}{2}\right)$$
for appropriate ladder operators $a,a^\dagger$ with $[a,a^\dagger] = 1$. The lowest eigenspace (lowest Landau level, corresponding to the eigenvalue $B$) can be represented as 
\begin{equation}\label{eq:LLL}
\mathrm{LLL} = \left\{ f(z) e^{-\frac{B}{4}|z|^2} \in L ^2 (\R^2), f \mbox{ holomorphic }\right\}
\end{equation}
and the $n$-th Landau level can be obtained as $\left(a^\dagger\right)^n \mathrm{LLL}.$ Hence each energy level is infinitely degenerate when working on the full plane. Well-known arguments indicate that this degeneracy is reduced in finite regions, with a degeneracy $\propto B \times \mbox{ Area }$. We give one such heuristic argument\footnote{Another one is that~\eqref{eq:Landau} can be restricted to a rectangle whose area is a multiple of $2\pi B^{-1}$, imposing magnetic-periodic boundary conditions see~\cite{AftSer-07,Almog-06,FouKac-11,Perice-22,NguRou-22} or~\cite[Sections~3.9 and 3.13]{Jain-07}. The energy levels are then the same as above, with degeneracy exactly $B (2\pi)^{-1} \times$ area of the rectangle.}. 

The orthogonal projector $\Pi_0$ on $\mathrm{LLL}$ can be expressed using vortex coherent states~\cite{ChaFlo-07,RouYng-19} in the form 
\begin{equation}\label{eq:Pi0}
 \Pi_0 = \frac{B}{2\pi} \int_{R\in \R^2} |\Psi_{0,R} \rangle \langle \Psi_{0,R}| d R 
\end{equation}
and similarly the orthogonal projector on the $n$-th Landau level is 
\begin{equation}\label{eq:Pin}
 \Pi_n =  \left(a^\dagger\right)^n \Pi_0 \left(a^\dagger\right)^n = \frac{B}{2\pi} \int_{R\in \R^2} |  \left(a^\dagger\right)^n \Psi_{0,R} \rangle \langle  \left(a^\dagger\right)^n \Psi_{0,R}| d R.
\end{equation}
The exact expression of $\Psi_{0,R}$ is known, but what matters here is that this function (as well as $\left(a^\dagger\right)^n \Psi_{0,R}$) is very localized (on the scale of the magnetic length $B^{-1/2}$) around the point $R\in\R^2$. Hence, to approximate the number of eigenstates of $H$ with eigenvalue $2B (n+1/2)$ localized in a given large domain $\Omega \subset \R^2$, it makes sense to restrict the integration in~\eqref{eq:Pi0}-\eqref{eq:Pin} to $\Omega$ and compute the trace of the so-obtained operator, namely 
$$ N_B (\Omega)= \frac{B}{2\pi} |\Omega|.$$      
This gives a good approximation of the rank of the operator $\1_{\Omega} \Pi_n H \Pi_n \1_{\Omega}$, hence of the number of orthogonal energy eigenstates with eigenvalue $\sim 2B (n+1/2)$ that one can fit in the domain $\Omega$.

\bigskip

\noindent\textbf{The integer quantum Hall effect.}~Some plateaux (left of Figure $1$) in $R_{xy}$/drops in $R_{xx}$ occur at integer values of $\nu$ and it is not surprising that something special should happen there (again, it is highly non-trivial to derive the specific signature of the ``something special''). This can be understood in a non-interacting electrons picture, taking only the Pauli exclusion principle into account. One assumes that the magnetic kinetic energy, proportional to $B$, is the main player and that all other energy scales in~\eqref{eq:start Hamil} are negligible against it. By this we mean that $W$ is dropped in~\eqref{eq:start Hamil} and that the only effect of $V$ is to essentially confine the gas to a domain $\Omega$.

As the name indicates, the filling factor measures the ratio of electron number to number of available one-body states in a given Landau level (see the above considerations, keeping in mind that $(2\pi)^{-1} h = c = e = 1$): 
$$ \nu = 2\pi \frac{\rho}{B} = 2\pi \frac{N}{|\Omega| B} \simeq \frac{N}{N_B (\Omega)}$$
if $N$ electrons are confined to the region $\Omega$ with density $\rho = N/|\Omega|$. In the ground state of an independent electron picture, one fills the eigenstates of~\eqref{eq:Landau} with one electron each, starting from the lowest one. At integer $\nu$, the $\nu$ lowest Landau levels are thus completely filled, and the others completely empty, a very rigid and non-degenerate situation. This rigidity is actually important in order to treat the energy scales other than $B$ perturbatively.  

\bigskip

\noindent\textbf{The fractional quantum Hall effect.}~Many plateaux however occur at particular \emph{rational} filling factors and are impossible to explain in an independent electrons picture. Laughlin's groundbreaking theory~\cite{Laughlin-83,Laughlin-87,Laughlin-99} explains why something special ought to occur at 
\begin{equation}\label{eq:Lau frac}
 \nu = \frac{1}{\ell}, \quad \ell \mbox{ an odd integer} 
\end{equation}
e.g. at the right-most plateau $\nu = 1/3$ of Figure~$1$, but also at $\nu = 1/5$, a fraction also observed in experiments ($\nu = 1/9$ and lower is not observed, while $\nu = 1/7$ is borderline). The $\nu = 1/3$ fraction is the first to have been observed~\cite{StoTsuGos-82}, and the most stable. Fractions from the principal Jain sequence (very prominent on the figure)
\begin{equation}\label{eq:Jain frac}
 \nu = \frac{p}{2p+1}, \quad p  \mbox{ an integer} 
\end{equation}
are explained in terms of the composite fermions theory~\cite{Jain-07}, a  generalization of Laughlin's theory we will not touch upon. Fractions of the form $$\nu = 1 - \frac{p}{2p+1}$$
are particle/hole symmetric pendants of the former~\eqref{eq:Jain frac}, and thus we cover most of the fractions seen on Figure~$1$. There are other, more exotic, fractions and features, but let us not get into that to focus on Laughlin's theory of the mother of all fractions, namely~\eqref{eq:Lau frac}.

\bigskip

\noindent\textbf{Restriction to the lowest Landau level}. We henceforth restrict to filling factors $\nu < 1$ . In the regime relevant to the quantum Hall effect, the gap $B$ between the magnetic kinetic energy levels is so large that the first approximation we make is to project all the physics down to as few Landau levels as possible. With filling ratio $\nu \leq 1$, the lowest Landau level is vast enough (again, see the above heuristics) to accommodate all particles, and thus we restrict available many-body wave-functions to those made entirely of lowest Landau\footnote{Generalizations to larger filling factors, when one works in an excited Landau level, are discussed in~\cite{RouYng-19}.} levels orbitals~\eqref{eq:LLL}. It is in fact convenient to work on the full space at first. The restrictions to finite area/density will actually be performed later, and we will have to make sure they are coherent with our aim: a thermodynamically large system with density $\rho \sim B \nu (2\pi)^{-1}$.

\bigskip

\noindent\textbf{Killing the interaction's singularity}. The main energy scale, the magnetic kinetic energy, is now frozen by projecting all one-body states to~\eqref{eq:LLL}. Laughlin's key idea is that the next energy scale to be considered is the pair interaction, and more precisely its singular short-range part. The wave-function~\eqref{eq:PsiLau} is introduced in order to reduce as much as possible the probability of particle encounters. Since we need the function to belong to
\begin{equation}\label{eq:LLLN}
 \mathrm{LLL}_N = \left\{ A(z_1,\ldots,z_N) e^{-\frac{B}{4} \sum_{j=1} ^N |z_j| ^2 }, \quad A \mbox{ analytic and antisymmetric} \right\}
\end{equation}
there is not much freedom. $\PsiLau$ is designed to vanish when $z_i=z_j$ while preserving the anti-symmetry and analyticity. It may seem that $\ell$ is a free variational parameter. But so far we thought somewhat grand-canonically: we have not fixed the density of our system yet. It turns out (this follows from Theorem~\ref{thm:plasma MF} below) that the one-particle density of Laughlin's function satisfies
\begin{equation}\label{eq:dens Lau}
 \varrho_{\PsiLau} (\bx) \simeq \frac{B}{2\pi \ell} \1_{|\bx| \leq \sqrt{\frac{2N\ell}{B}}}.
\end{equation}
That is, it lives on a thermodynamically large length scale\footnote{That the support is a disk is of no importance, since bulk properties of systems in thermodynamic limits are typically independent of the sample's shape.} and has filling factor $\nu = \ell^{-1}$ (if you bear with me concerning the choice of units in~\eqref{eq:mag hamil}). 

Now we can answer our original question ``what is special about filling factor $\nu = \ell ^{-1}$ ?'' The answer is that, at such parameter values, we may form a Laughlin state of exponent $\ell$ as approximate ground state of our system. It minimizes the magnetic kinetic energy exactly, and does a very good job at reducing the short-range part of the interaction.

\bigskip

\noindent\textbf{Towards a rigorous derivation of Laughlin's function.} The second part of the above derivation is very heuristic, and will probably stay that way in the case of the true 3D Coulomb interaction. However, if one is willing to approximate the short-range part of the interaction as a sharply peaked delta-like potential, one may indeed derive rigorously the Laughlin state (and/or variants) in a physically relevant limit. This is based on the fact that the Laughlin function is an \emph{exact} ground state for an approximate interaction of zero-range, projected on the lowest Landau level. For a precise formulation of this, and the derivation of such model interactions from scaled ones, I refer to~\cite{LewSei-09,SeiYng-20}. A Gross-Pitaevskii-like limit of bosonic models based on effective delta interactions projected in the lowest Landau level is studied in a very nice paper by Elliott and co-workers~\cite{LieSeiYng-09}.

The main open problem in this direction is to make the derivation of Laughlin's function alluded to above \emph{uniform in the particle number} $N$. This depends on a \emph{spectral gap conjecture} for effective zero-range interactions, whose formulation can be found in~\cite[Appendix]{Rougerie-xedp19} and references therein. Partial progress towards the conjecture are in~\cite{NacWarYou-20a,NacWarYou-20b,WarYou-21,WarYou-21b}.

\bigskip

\noindent\textbf{Laughlin quasi-holes}. So far we have argued that Laughlin's function is a good ansatz for the ground state of the system at the relevant filling factor, when neglecting the effect of the external potential $V$ and the long-range part of the interaction $W$ in~\eqref{eq:start Hamil}. That is not the end of the story, for the latter ingredients do exist in actual experiments, in particular, the disorder landscape that impurities enforce in $V$ is crucial to the quantum Hall effect.

The Laughlin state should in fact be seen as the ``vacuum'' of a theory explaining the FQHE experimental data. The next step is to construct the quasi-particles generated from said vacuum when suitably moderate external fields are applied, such as those generating the currents in experiments. 

It is in fact easier to argue about quasi-holes, generated e.g. when the filling factor is lowered a little from the magic fraction $\ell ^{-1}$, as when moving towards the right on Figure~$1$. The salient feature is that we stay on the same FQHE plateau for a while when doing so. It must hence be that the ground state of the system stays ``Laughlin-like'' for reasonably smaller $\nu$. In fact, Laughlin's next key idea is two-fold 
\begin{itemize}
 \item for smaller filling factors, the ground state is generated from~\eqref{eq:PsiF} by adding uncorrelated quasi-holes as in~\eqref{eq:Psif}-\eqref{eq:polynomial}. These are typically pinned by the impurities of the sample (modeled by $V$ in~\eqref{eq:start Hamil}).
 \item when applying an external field at $\nu$ close to $\ell^{-1}$, the current is carried by the motion of such quasi-holes. 
\end{itemize}
The second idea in particular is quite far-reaching: it has by now been measured~\cite{SamGlaJinEti-97,YacobiEtal-04,MahaluEtal-97} that the current is carried in fractional lumps of $e \ell^{-1}$ and~\cite{BarEtalFev-20,NakEtalMan-20} that the charge carriers obey fractional quantum statistics, i.e. are emergent anyons~\cite{AroSchWil-84,Halperin-84,LunRou-16,LamLunRou-22,ZhaSreGemJai-14,CooSim-15}.

\bigskip

\noindent\textbf{Stability of the Laughlin phase}. The last point motivates the variational problem studied in Section~\ref{sec:prob}. The model incorporates the ingredients from~\eqref{eq:mag hamil} that are not frozen by the aforementioned reductions: the external potential $V$ representing trapping and disorder and the long-range part of the interaction potential $W$. By restricting the variational set as in~\eqref{eq:qm_energy} we take for granted the basic ingredients sketched above, but back-up a little by noting that they actually point to the general form~\eqref{eq:PsiF} for trial states. The Laughlin function~\eqref{eq:PsiLau} and associated quasi-holes states~\eqref{eq:Psif} certainly are the simplest, and hence the first to try in order to explain experimental data. But ideally they should be singled out from the full set~\eqref{eq:PsiF} by minimizing\footnote{This is like in degenerate perturbation theory: we minimize the smallest energy scales~\eqref{eq:start Hamil} amongst all possible ``minimizers''~\eqref{eq:PsiF} of the largest energy scales.} the remaining energy scales in~\eqref{eq:many_body_energy}. This is what Theorem~\ref{thm:ener} proves, under some simplifying assumptions that we now discuss.

In scaling the external potential as in~\eqref{eq:rescaled_V} we make it live on the natural, thermodynamically large, length-scale of the Laughlin function. This is very reasonable for the trapping part of the potential, but much less so for the part modeling disorder, which typically lives on a much shorter length scale. In fact the shortest length scale we could allow is that dictated by Theorem~\ref{thm:incomp}, so it does not need to be thermodynamically large. Improving it to realistic values however remains an open problem, and we prefer for simplicity to work on a single length scale in order not to obscure the main statements. 

In Theorem~\ref{thm:ener} we assume the interaction to be smooth. This is because it is supposed to represent the long-range part only, the singular short-range part being taken care of by restricting to~\eqref{eq:PsiF}. Scaling $W$ as in~\eqref{eq:rescaled_w} has the merit of making the two terms in~\eqref{eq:many_body_energy} of the same order of magnitude, as in a mean-field limit. This also simplifies statements a lot, but for interactions scaling like 3D Coulomb, this is actually the correct thing to do, see~\cite[Section~2.2]{OlgRou-19}.

Concerning the smallness assumption on $\lambda$ in~Theorem~\ref{thm:ener}, it corresponds to the fact that the filling factor should stay close to $\ell^{-1}$ for the theorem to be true. Too large a deviation makes the system jump to a different FQHE plateau, e.g. a Laughlin state with higher exponent. What is slightly tricky is that we do not work at fixed density but fixed particle number. But increasing the (repulsive) interaction strength has the net effect of spreading the system further, and hence lowering the density (see again~\cite[Section~2.2]{OlgRou-19} for more details). An upper bound on $|\lambda|$ is thus necessary for the statement to hold. We do not however provide a meaningful estimate of the size of $|\lambda|$ needed for the proof to carry through, which is probably model-dependent.

\section{The one-component plasma}\label{sec:jellium}

Our ultimate goal is to get to grips with the many-body density $|\Psi_F|^2$ of functions defined as in~\eqref{eq:PsiF}. Since this will be made possible by the plasma analogy~\eqref{eq:plasma}-\eqref{eq:class hamil}, we first briefly review the statistical mechanics of classical Coulomb systems (referring to~\cite{Serfaty-14,Serfaty-15,Serfaty-17,Lewin-22} for more complete accounts). The challenging step of including a general many-body analytic factor $F$ as in~\eqref{eq:class hamil} will mostly be dealt with in the next section. We first focus on a case closer to the target $F = f^{\otimes N}$ appropriate to describe quasi-holes wave-functions~\eqref{eq:Psif}. 

In view of~\eqref{eq:crazy pot}-\eqref{eq:qhonebod}, a general Hamiltonian including~\eqref{eq:class hamil} with $F = f ^{\otimes N}$ as a particular case is as follows:
\begin{equation}\label{eq:Coulomb system}
H_N (\bx_1,\ldots,\bx_N) := \sum_{j= 1} ^N V (\bx_j) + \sum_{1\leq j < k \leq N} \coul (\bx_j-\bx_k) 
\end{equation}
where $\bx_1,\ldots,\bx_N \in \R^d$ are coordinates of particles in the Euclidean space, $V:\R^d \mapsto \R$ is an external potential and 
\begin{equation}\label{eq:kernel}
\coul (\bx) = \begin{cases}
               \frac{1}{|\bx| ^{d-2}} \mbox{ if } d\geq 2\\
               -\log |\bx| \mbox{ if } d = 2.
              \end{cases}
\end{equation}
We only consider space dimensions $d\geq 2$ in the sequel, where $\coul$ is the fundamental solution of Laplace's equation:
\begin{equation}\label{eq:Laplace}
 -\Delta \coul = |\mathbb{S}^{d-1}| \delta_0
\end{equation}
so that the potential $\varphi_\rho$ generated by a charge distribution $\rho$ is obtained through 
\begin{equation}\label{eq:Laplace2}
 \varphi_\rho  (\bx) = \int_{\R^d} \coul(\bx-\by) \rho (\by) d\by.
\end{equation}
We will be interested in equilibrium properties, namely in the Gibbs states at temperature $T>0$
\begin{equation}\label{eq:Gibbs}
\mub_{T,N} := \frac{1}{\cZ (T,N)} \exp\left( -\frac{1}{T} H_N \right)
\end{equation}
minimizing the free-energy
\begin{equation}\label{eq:free ener func}
 \cF_{T,N} [\nu] := \int_{\R^{dN}} \left( H_N \nub + T \nub \log \nub \right)  
\end{equation}
amongst probability measures $\nub$ over $\R^{dN}$. As usual the infimum is given in terms of the partition function $\cZ_{T,N}$ normalizing~\eqref{eq:Gibbs} as 
$$ F (T,N) = - T \log \cZ (T,N)$$
and we identify the $T=0$ problem with the ground state. Namely 
$$ F (0,N) = \inf_{\R^{dN}} H_N$$
and $\mub_{0,N}$ is the empirical measure associated with a minimum point.

\subsection{Homogeneous systems and the thermodynamic limit}\label{sec:thermo} 

We start by describing a fundamental contribution by Elliott Lieb and Heide Narnhofer~\cite{LieNar-76}, inspired by the methods of~\cite{LebLie-69,LieLeb-72}. Namely we confine the Coulomb gas described above to a finite container of volume $L^d$, fix the density
$$\rho := \frac{N}{L^d}$$
and take the thermodynamic limit $N,L\to \infty$. For this to make sense, we need to make the system neutral. In the jellium model this is done by taking an external potential generated by a constant neutralizing background of  density $-\rho$. Hence we think of one species of charges as fixed and spread in space, and the other as point-like and moving in the ``jelly'' thus generated. 

The highly non-trivial point is to quantify the screening of the background by the mobile point charges, leading to a system neutral on length scales much larger than the microscopic typical inter-particle distance $N^{-1/d}$. Hence the long-range tail of the Coulomb potential is not felt all across the box and the free energy can be extensive (after, of course, having taken the energy of the background into account), as indicated by the 

\begin{theorem}[\textbf{Thermodynamic limit for jellium}]\mbox{}\label{thm:thermo}\\
Let $\Omega$ be a regular simply connected set and $\Omega_L$ its dilation by a factor $L$. Let $\rho \in \R^+$, $V=V_L$ be given by 
\begin{equation}\label{eq:pot hom}
V_L (\bx):= 
\begin{cases} \displaystyle
- \rho \int_{\Omega_L} \coul (\bx-\by)d\by \mbox{ if } \bx \in \Omega_L\\
+ \infty \mbox{ otherwise}
   \end{cases}
\end{equation}
and the associated energy by
$$ 
E (L,\rho) := \frac{\rho^2}{2} \iint_{\Omega_L \times \Omega_L} \coul (\bx-\by)d\bx d\by.
$$
The limit 
\begin{equation}\label{eq:therm lim}
f (T,\rho) := \displaystyle\lim_{L,N \to \infty, \frac{N}{L^d} \to \rho} \frac{\displaystyle F (T,N) + E (L,\rho) }{L^d} 
\end{equation}
exists and is independent of the shape of $\Omega$.
\end{theorem}

This follows~\cite{LieNar-76,SarMer-76} from the method of~\cite{LebLie-69,LieLeb-72} which uses Newton's theorem and averages over rotations around well-chosen centers to quantify screening. A difficulty is that the background is fixed, so that extra care has to be taken with this procedure compared to~\cite{LebLie-69,LieLeb-72} where the background is replaced by point charges. On the other hand, stability of matter is not an issue in this set-up, so that one can deal with the classical model. In~\cite{LebLie-69,LieLeb-72} one has to use the Heisenberg and Pauli principles of quantum mechanics in a highly non-trivial way~\cite{Lieb-76,LieSei-09} to prove that even the $\liminf$ makes sense. 

\subsection{Inhomogeneous systems and the mean-field approximation}\label{sec:MF}

A natural follow-up question is to choose a inhomogeneous background distribution of charge in~\eqref{eq:Coulomb system}. Thus we now choose a general external potential $V$. It is still desirable that the system lives on a thermodynamic length scale $\sim N^{1/d}$, and that the energy be extensive. We achieve this by picking 
\begin{equation}\label{eq:inhom}
 V (\bx) = V_N (\bx) = N^{\frac{2}{d}} v (N^{-1/d} \bx)
\end{equation}
with $v$ a fixed confining potential (i.e growing at infinity). Then, changing length units by setting $\bx_j = N^{1/d} \by_j$ in~\eqref{eq:Coulomb system} we obtain an energy in mean-field scaling\footnote{\label{foo:delta}Hereafter $ \delta_{d=2}= 1$ in 2D and $0$ otherwise.}
$$ H_N (\by_1,\ldots,\by_N) = N^{\frac{2}{d}} \left( \sum_{j= 1}^N v(\by_j) + \frac{1}{N} \sum_{1\leq j < k \leq N} \coul(\by_j-\by_k)\right) - \frac{N(N-1)\log N}{4} \delta_{d=2}.$$
In the rescaled Hamiltonian in parenthesis, the two terms formally weigh the same and thus the $\by_j$'s will want to stay in a domain of fixed volume. The latter, and the overall shape of the density is obtained by minimizing a continuum/uncorrelated version of the above:
\begin{equation}\label{eq:MF func}
\cEMF [\sigma] = \int_{\R^d} V \sigma (\bx) d\bx + \frac{1}{2} \iint_{\R^d \times \R^d} \sigma (\bx) \coul(\bx-\by)\sigma (\by) d\bx  d\by.
\end{equation}
In fact, rescaling lengths in~\eqref{eq:Gibbs}-\eqref{eq:free ener func}, one sees that the temperature is effectively small, so that the entropy does not appear in the leading order of the energy (it does~\cite{MesSpo-82,Kiessling-89,Kiessling-93,Kiessling-09b,CagLioMarPul-92,CagLioMarPul-95,Serfaty-20,ArmSer-20,ArmSer-21} if one allows for $T$ to grow appropriately fast when $N\to \infty$, see also~\cite[Chapter~2]{Rougerie-cdf,Rougerie-LMU}). 

\begin{theorem}[\textbf{Mean-field limit for inhomogeneous Coulomb systems}]\label{thm:plasma MF}\mbox{}\\
In the set-up just described, with $v$ a sufficiently (requirements are low) regular function growing polynomially (not really required) at infinity, and $T\geq 0$ fixed, we have that (cf Footnote~\ref{foo:delta}) 
\begin{equation}\label{eq:MF lim ener}
F(N,T) = -\frac{N(N-1)(\log N)}{4}\delta_{d=2} + N^{1+2/d} \EMF (1+o(1)) 
\end{equation}
where $\EMF$ is the infimum of~\eqref{eq:MF func} amongst probability measures on $\R^d$. 

Moreover, let $\varrho_N^{(k)}$ be the $k$-particle density of the full Coulomb system\footnote{In particular $\varrho_N^{(1)} (\bx) = \sum_{j=1}^N \delta_{\bx = \bx_j}$ with $\bx_1,\ldots,\bx_N$ a minimum point for $H_N$ if $T=0$.}, i.e.
$$ \varrho_N^{(k)} (\bx_1,\ldots,\bx_k) = \frac{N!}{(N-k)!} \int_{\R^{d(N-k)}} \mub_{T,N}(\bx_1,\ldots,\bx_N) d\bx_{k+1} \ldots d\bx_N$$
with $\mub_{T,N}$ the Gibbs measure~\eqref{eq:Gibbs} if $T>0$.  Then 
\begin{equation}\label{eq:MF lim dens}
 \varrho_N^{(k)} (N^{1/d} \by_1,\ldots,N^{1/d}\by_k) \underset{N\to \infty}{\rightharpoonup} \rhoMF (\by_1) \ldots  \rhoMF (\by_k)
\end{equation}
weakly as measures, where $\rhoMF$ is the unique (using that $\coul(\,.\,)$ is of positive type) minimizer for $\cEMF$.   
\end{theorem}

We did not aim at the greater generality or precision in the above. At various degrees of both these criteria, proofs may be found in~\cite{AndGuiZei-10,Forrester-10,Serfaty-15,RouSerYng-13b,RouSer-14,SafTot-97,BenZei-98,BenGui-97,Hardy-12,ChaGozZit-13,ChaHarMai-16} and many related sources. This vindicates~\eqref{eq:dens Lau}, for in this case it is particularly easy to compute $\rhoMF$. 

To see more precisely why the entropy does not contribute at this order, observe that one can expect the mean-field ``uncorrelated'' behavior   
\begin{equation}\label{eq:ansatz}
\mub_{T,N} (\bx_1,\ldots,\bx_N) \simeq N^{-N}  \rhoMF \left(N^{-1/d}\bx_1\right) \ldots \left(N^{-1/d}\bx_N\right). 
\end{equation}
In fact, the above result is compatible with this ansatz, and shows there is a lot of truth in it. The entropy in such an ansatz is 
\begin{equation}\label{eq:small ent}
- \int_{\R^{dN}} \mub_{T,N} \log \mub_{T,N} \simeq  N \log N - N \int_{\R^d} \rhoMF \log \rhoMF.
\end{equation}
Since the temperature $T$ in~\eqref{eq:free ener func} is fixed, the contribution of $-T \times$ entropy to the free energy is much smaller than the terms identified in~\eqref{eq:MF lim ener}.

\subsection{Inhomogeneous systems and the local density approximation}\label{sec:LDA}

Informally, Theorem~\ref{thm:thermo} is a very precise version of Theorem~\ref{thm:plasma MF} in the particular case (homogeneous system) to which it applies: the next-to-leading order beyond mean-field is identified. The equivalent result for Inhomogeneous systems was obtained much later, in~\cite{SanSer-14,RouSer-14,PetSer-14,RotSer-14,PetRot-16} at $T=0$ and in~\cite{ArmSer-20,LebSer-15} at $T>0$. Later developments can be found e.g. in~\cite{Serfaty-20,ArmSer-20,ArmSer-21,Leble-16,Leble-15b,LebSer-16,BauBouNikYau-15,BauBouNikYau-16}. The formulation of the result is in a somewhat different spirit from Theorem~\ref{thm:thermo} in these references (and many more things are proved beyond what we state), but we refer to~\cite{Lewin-22} for an explanation of the fact that, indeed, the statement below follows from~\cite{SanSer-14,RouSer-14,PetSer-14,LebSer-15}:

\begin{theorem}[\textbf{Local density approximation for inhomogeneous Coulomb systems}]\label{thm:LDA}\mbox{}\\
Under the same assumptions as in the previous theorem, we have, in the limit $N\to \infty$ with $T$ fixed 
\begin{equation}\label{eq:LDA}
F(N,T) = - \frac{N^2\log N}{4}\delta_{d=2} + N^{1+2/d} \EMF  + N \int_{\R^d} f\left(T,\rhoMF(\bx)\right) d\bx + o (N)
\end{equation}
where $f(T,\rho)$ is defined in Theorem~\ref{thm:thermo} and $ \delta_{d=2}= 1$ in 2D and $0$ otherwise.
\end{theorem}

This is called a ``local density approximation'' because the correction to mean-field theory is obtained by integrating the free-energy density of the homogeneous system at density $\rhoMF(\bx)$ over $\bx$. This means that, locally at the microscopic scale, the system is in thermal equilibrium at the density set by the macroscopic mean-field theory. This separation of scales is again a powerful manifestation of screening in Coulombic matter. See~\cite{LewLieSei-18,LewLieSei-19a,LewLieSei-19b,CotPet-17,CotPet-19} for similar results in the context of the uniform electron gas and density functional theory. It is noteworthy that 
\begin{itemize}
 \item the fixed temperature shows up at the level of precision of the above but is absent from the leading order in Theorem~\ref{thm:plasma MF}. 
\item precise estimates on the remainder are obtained in~\cite{ArmSer-20}. The order of magnitude thereof are presumably optimal, for they scale precisely like boundary terms (at least in the homogeneous case where Theorem~\ref{thm:LDA} reduces to Theorem~\ref{thm:thermo}). 
\end{itemize}

\subsection{Renormalized Jellium energy}

In the proof of Theorem~\ref{thm:LDA} (and for the derivation of important corollaries not mentioned here), it is useful to characterize the homogeneous Jellium's free energy not as the \emph{thermodynamic limit of an infimum in finite volume}, but as the \emph{infimum of a quantity directly defined in infinite volume}. We briefly sketch this below, referring to~\cite{SanSer-12,SanSer-14,RouSer-14,PetSer-14,LebSer-15} for more details. That the quantities defined below coincide with those of Section~\ref{sec:thermo} follows from the fact that the results of~\cite{SanSer-12,SanSer-14,RouSer-14,PetSer-14,LebSer-15,ArmSer-20}, bearing firstly on the inhomogeneous setting of Sections~\ref{sec:MF}~and~\ref{sec:LDA}, apply as well in the homogeneous setting by choosing the external potential as in Theorem~\ref{thm:thermo} (see~\cite{Lewin-22}, in particular Remark~38 therein for more comments on this point). 

A first key point is to define the energy of a charge configuration via the electric field it generates. 

\begin{defi}[\textbf{Admissible electric fields}]\label{def:adm field}\mbox{}\\ 
Let $\rho>0$. Let $\j$ be a vector field in $\R^d$. We say that $\j$ belongs to the class $\bam$ if $\j = \nabla h$ with
\begin{equation}\label{curlj}
-\Delta h = c_d \Big(\sum_{p\in \Lambda}N_p \delta_p - \rho \Big) \quad \text{in} \ \R^d
\end{equation}
for some discrete set $\Lambda \subset \R^d$, and $N_p$ integers in $\N^*$.
\end{defi}

One should think of $\j$ as the electric field (i.e. gradient of the potential) generated via Laplace's equation~\eqref{eq:Laplace} by a configuration of point charges and a uniform background of density $\rho$. The ``renormalization'' alluded to in this subsection's title enters via the smearing of point charges on a length scale $\eta$, ultimately sent to $0$ after the subtraction of appropriate counter-terms.

\begin{defi}[\textbf{Smeared electric fields}]\label{def:smear field}\mbox{}\\
Pick some arbitrary fixed {\it radial} nonnegative function $\sigma$, supported in $B(0,1)$ and with integral $1$. For any point $p$ and $\eta>0$ we introduce the smeared charge 
\begin{equation}\label{eq:smeared charge}
\delta_p^{(\eta)}= \frac{1}{\eta^d}\sigma \left(\frac{x}{\eta}\right) *  \delta_p.
\end{equation}
Let the radial function $g_\eta$ be the unique solution (Newton's theorem~\cite[Theorem~9.7]{LieLos-01} is used here) to
\begin{equation}\label{eqf0}
\begin{cases}
 -  \Delta g_\eta= c_d\left( \delta_0^{(\eta)} - \delta_0 \right) \quad \text{in} \ \R^d
\\   
 g_\eta\equiv 0 \quad  \text{in} \ \R^d \backslash B(0,\eta).
\end{cases}
\end{equation}
For any vector field $\j=\nabla h$ satisfying
\begin{equation}\label{dij}
-\mathrm{div}\, \j= c_d\Big(\sum_{p \in \Lambda}N_p \delta_p -\rho\Big)
\end{equation}
in some subset $U$ of $\R^d$, with $\Lambda \subset U$ a discrete set of points, we let
$$\j_\eta := \nabla h + \sum_{p\in \Lambda} N_p \nabla g_\eta(x-p) \qquad h_\eta= h + \sum_{p \in \Lambda} N_p g_\eta(x-p).$$
We have
\begin{equation}
\label{delp}
-\mathrm{div}\,  \j_\eta  = - \Delta h_\eta =  c_d\Big(\sum_{p \in \Lambda}N_p \delta_p^{(\eta)} -\rho\Big).\end{equation}
\end{defi}

To define the energy (per volume) of an infinite configuration of point charges minus uniform background, we observe that formally (that is, modulo subtracting the infinite self-energies of the point charges) it ought to be given (using~\eqref{eq:Laplace} again) by the mean value
$$ \dashint_{\R^d} |\j| ^2$$
where $\j$ is the electric field (see Definition~\ref{def:adm field}). This is where the renormalization via screening takes place: we smear charges, consider 
$$ \dashint_{\R^d} |\j_\eta| ^2$$
as in Definition~\ref{def:smear field}, remove the self-energies of individual smeared charges, and then pass to the limit $\eta \to 0$. A key point in the following definition is that we pass to the infinite volume limit \emph{before} letting $\eta \to 0$.

\begin{defi}[\textbf{The renormalized jellium energy}]\label{def:renorm ener}\mbox{}\\
Let $\sigma$ be as in Definition~\ref{def:smear field} and the energy of charges smeared on radius $1$ be
\begin{align}\label{kapd}
\kappa_d&= c_d \iint_{\R^d \times \R^d} \sigma (\bx) \coul(\bx-\by) \sigma(\by) d\bx d\by \quad \text{for} \  d\ge 3,  \quad \kappa_2=c_2\nonumber\\ 
\qquad \gamma_2&=  c_2 \iint_{\R^d \times \R^d} \sigma (\bx) \coul(\bx-\by) \sigma (\by) d\bx d\by  \ \text{for} \ d=2.
\end{align} 
with
\begin{equation}\label{defc}
c_2 = 2\pi, \qquad c_d = (d-2)|\mathbb{S}^{d-1}| \ \text{when} \ d\ge 3
\end{equation}
For any $\j \in  \bam$, we define ($K_R$ is a hypercube of side-length $R$, $\dashint_{K_R} F$  is the mean of $F$ over it)
\begin{equation}\label{We}
\W_\eta(\j) = \limsup_{R\to \infty} \dashint_{K_R}  |\j_\eta|^2 - \rho \left(\kappa_d  \coul(\eta)+\gamma_2 \delta_{d=2}\right)
\end{equation}
and the renormalized jellium energy is given by 
\begin{equation}\label{weta}
\W(\j) = \liminf_{\eta\to 0}\W_\eta(\j) =\liminf_{\eta\to 0}\left(\limsup_{R\to \infty} \dashint_{K_R}   |\j_\eta|^2 - \rho \left(\kappa_d \coul(\eta)+\gamma_2 \delta_{d=2}\right)\right).
\end{equation} 
Again, $ \delta_{d=2}= 1$ in 2D and $0$ otherwise.
\end{defi}

A first statement we can make is that 
\begin{equation}\label{eq:infvol T0}
\boxed{\inf_{\j \in \bam} \W (\j) = f(0,\rho)}
\end{equation}
where $f(0,\rho)$ is the ground-state energy per volume defined by Theorem~\ref{thm:thermo}. This follows from the aforementioned works~\cite{SanSer-14,RouSer-14,PetSer-14,RotSer-14,PetRot-16}  by choosing the external potential as in Theorem~\ref{thm:thermo} (see also~\cite{Lewin-22}). An extension to positive temperatures requires more definitions~\cite{LebSer-15}, for which we shall be somewhat less precise. We call a point process $P$ a probability measure over locally finite point configurations in $\R^d$, or equivalently over the set of non-negative, purely atomic Radon measures on $\R^d$ giving an integer mass to singletons.

Then we have 

\begin{definition}[\textbf{Energies and entropies of point processes}]\label{defi:jellium free ener}
Let $\Lambda$ be a point configuration. Its renormalized jellium energy with background $\rho>0$ is\footnote{The integers $N_p$ in~\eqref{curlj} account for the fact that there might be multiple points in the configuration $\Lambda$ seen as a set.}   
$$ \mathbb{W} (\Lambda,\rho) := \frac{1}{c_d}\inf \left\{ \W(\j), \j \in \bam \mbox{ satisfying~\eqref{curlj} for that particular $\Lambda$}\right\}$$
Let $P$ be a point process. Its renormalized jellium energy with background $\rho>0$ is
$$ \widetilde{\mathbb{W}} (P,\rho) = \left\langle \mathbb{W} (.,\rho) \right\rangle_P$$
with $\left\langle \, .\, \right\rangle_P$ denoting expectation in $P$.

Let $\Pi^m$ be the Poisson point process with intensity $\rho>0$. For a stationary (translation-invariant) point process $P$ we let its entropy relative to the Poisson point process be
$$ \mathrm{Ent} (P,\Pi^\rho) := \lim_{R\to \infty} \frac{1}{R^d} \mathrm{Ent} \left(P_{K_R}, \Pi^\rho_{K_R}\right)$$
where $P_{K_R}, \Pi^\rho_{K_R}$ are the restriction of the process to the hypercube $K_R$, and $\mathrm{Ent}$ is the usual relative entropy of probability measures defined over the same probability space. I.e, in terms of the Radon-Nikodym of $\mu$ with respect to $\nu$ 
$$ \mathrm{Ent} (\mu,\nu) = \int \frac{d\mu}{d\nu}\log \left( \frac{d\mu}{d\nu}\right) d\nu $$
with the convention is $+\infty$ is the Radon-Nikodym derivative does not exist.
\end{definition}

The positive $T$ equivalent of~\eqref{eq:infvol T0} is (as it should)
\begin{equation}\label{eq:infvol T}
\boxed{\inf \left\{ \widetilde{\mathbb{W}} (P,\rho) + T \mathrm{Ent} (P,\Pi^\rho), P \mbox{ a point process over } \R^d\right\} = f(T,\rho)}
\end{equation}
where $f(T,\rho)$ is the free-energy per volume defined by Theorem~\ref{thm:thermo}. The identity follows from~\cite{LebSer-15,ArmSer-20} in the same way as~\eqref{eq:infvol T0}. There would be a lot more to say about these definitions (including the extension to tagged point processes, crucial to~\cite{LebSer-15}) and their uses. To keep things between bounds, we stick to the following comments.

The infinite volume quantities defined above are very useful to obtain estimates and limit theorems on the fluctuations around the mean-field Theorem~\ref{thm:plasma MF}: central limit theorems, large deviation principles ... derived in the aforementioned references. The formulation via the electric field in Definition~\ref{def:renorm ener} permits to quantify screening mechanisms differently from what was alluded to above, i.e. without appealing to local rotational invariance and then Newton's theorem. Briefly, ``good'' configurations (energy-wise) of point charges can be modified slightly in order not to change the energy or entropy too much, while making the associated electric field vanish outside of a large hyper-rectangle. Modified configurations in neighboring hyper rectangles can then be glued together without introducing divergences in the field at the interface, and thus making the energies add up. Thereby one obtains trial states for large domains by gluing equilibrium configurations in smaller domains, which is key to a form of additivity allowing to deduce the existence of thermodynamic quantities.

\subsection{Separation of points in the ground state}

There is one (unpublished) result of Elliott's which, in addition to its independent interest, played an important role in putting the tools of the previous section to good use. To see this, observe that it is not obvious from Definition~\ref{def:renorm ener} that the infimum of the infinite volume renormalized jellium energy~\eqref{weta} is finite. We have to make sure that the negative, diverging when $\eta\to 0$, counter-terms do indeed cancel corresponding infinities in the main term. In the approach of~\cite{RouSer-14} (simplifying~\cite{SanSer-14}, and in turn simplified in~\cite{PetSer-14} and~\cite{ArmSer-20} using other ingredients) this is achieved as follows:
\begin{itemize}
 \item for a lower bound it is sufficient to consider (quasi-)minimizers.
 \item by a variant of Elliott's argument, point charges corresponding to (quasi-)minimizers via~\eqref{curlj} are well-separated.
 \item if points are more than a distance of order $\eta$ apart, the \emph{radial} smearing of point charges in~Definition~\ref{def:smear field} does not change the interaction between different points, by Newton's theorem.  
 \item a few calculations and estimates then vindicate that indeed infinities compensate one-another in the limit $\eta \to 0$.
\end{itemize}
The crucial ingredient, in the second point above, shows that the minimal distance between points from quasi-minimizing configurations  is bounded below, uniformly in the limit $R\to \infty$ (thermodynamic limit) of~\eqref{We} and in the limit $\eta \to 0$. This is obtained by a variant (mostly accounting for the positive smearing parameter $\eta$) of the following (see also~\cite[Lemma~24]{Lewin-22}).

\begin{theorem}[\textbf{Separation of charges} (Lieb, unpublished)]\label{thm:sep}\mbox{}\\
 Same setting as in Theorem~\ref{thm:thermo}, with temperature $T=0$ and density 
 $$ \rho = \frac{N}{L^d}.$$
Let $(\bx_1,\ldots,\bx_N) \in \Omega_L^N$ be a minimizing configuration for the energy $F(0,N)$, pick a point therefrom, denoted $\bx_0$ without loss of generality (the energy is invariant under change of labels). Let $\delta$ be the radius of the ball of unit volume in $\R^d$. if 
$$ \mathrm{dist}\left(\bx_0,\partial \Omega_L\right) \geq \delta \rho ^{-1/d}$$
then 
$$ \min_{j\neq 0}\left|\bx_0 - \bx_j\right| \geq \delta \rho ^{-1/d}.$$
\end{theorem}

Even though the use of this theorem has now been by-passed~\cite{ArmSer-20,PetSer-14} to bound energies of the type of Definition~\ref{def:renorm ener} from below, variants of it are still crucial to prove separation/equidistribution of charge in related systems~\cite{PetSer-14,RotSer-14,PetRot-16}.

The proof is as short as it is elegant.

\begin{proof}
There is a ball of center $\bx_0$ and radius $\delta \rho ^{-1/d}$ fully included in $\Omega_L$. Assume for contradiction that there is at least another point $\bx_k\neq \bx_0$ in said ball, and consider variations of the energy with respect to the motion of that point, all the others being fixed. For the configuration to be a minimizer, it must be that $\bx_k$ sits at a minimum of the potential 
$$\Phi (\bx) := \sum_{j\neq k} \coul (\bx_k - \bx_j) - V_L (\bx)$$  
generated by the other points and the background, with $V_L$ as in~\eqref{eq:pot hom}. We split $\Phi$ as 
\begin{align}\label{eq:Lieb}
 \Phi (\bx) &= \Phi_1 (\bx) + \Phi_2 (\bx)\nonumber\\
 &= \coul (\bx - \bx_0) - \rho \int_{B(0,\delta \rho ^{-1/d})} \coul(\bx-\by)d\by \nonumber\\
 &+ \sum_{j\neq k,0} \coul (\bx - \bx_j) - \rho \int_{\Omega_L \setminus B(0,\delta \rho ^{-1/d})} \coul(\bx-\by)d\by.
\end{align}
By virtue of~\eqref{eq:Laplace}-\eqref{eq:Laplace2}, the potential $\Phi_2$ on the third line is superharmonic on $B(0,\delta \rho ^{-1/d})$. Indeed, the terms in the sum are superharmonic everywhere, and the other term is harmonic on $B(0,\delta \rho ^{-1/d})$. Hence, by the maximum principle, 
$$\Phi_2 (\bx) \geq \min_{\bx \in \partial B(0,\delta \rho ^{-1/d})} \Phi_2 (\bx)$$ 
for all $\bx \in B(0,\delta \rho ^{-1/d})$. On the other hand, the charge generating $\Phi_1$ on the second line is radial around $\bx_0$. Hence $\Phi_1$ has the same symmetry, and can be computed explicitly, by Newton's theorem~\cite[Theorem~9.7]{LieLos-01}. With $\delta$ chosen as in the statement, one can see that it takes its minimum (namely, $0$) on $\partial B(0,\delta \rho ^{-1/d})$ and is strictly positive in the interior of the ball. Hence there must exist a point $\by$ on $\partial B(0,\delta \rho ^{-1/d})$ with 
$$ \Phi (\by) < \Phi (\bx) \mbox{ for all } \bx \mbox{ in the interior of } B(0,\delta \rho ^{-1/d}).$$
$\Phi$ can be decreased by moving $\bx_k$ to $\by$, which contradicts the fact that our configuration was assumed to be a minimizer.
\end{proof}

\section{Incompressibility bounds for Coulomb ground states}

We turn to explaining how one can extend the idea of the proof of Theorem~\ref{thm:sep} to prove density bounds of the form~\eqref{eq:incomp vague}. Using the plasma analogy described in Section~\ref{sec:analogy}, this translates to density bounds on Gibbs equilibria of generalized 2D classical Coulomb systems. It turns out that the effective temperature in the plasma analogy is quite small in the limit $N\to \infty$: in fact one is exactly in the scaling described in Sections~\ref{sec:MF}-\ref{sec:LDA}. Our proof of Theorems~\ref{thm:incomp}-\ref{thm:exponential_bound} proceeds from bounds for the ground state at $T=0$, coupled with rough estimates relating Gibbs to ground states for small $T$. The latter part is that where we have to restrict to the non-optimal length scales in our statements, and for which an improvement would be most desirable. But this remains an open problem.  

In this note we restrict to explaining how to obtain density upper bounds for ground states of 2D classical Coulomb Hamiltonians of the form~\eqref{eq:class hamil}. Modulo a change of length and energy units, we can consider the following general Hamiltonian 
\begin{equation}\label{eq:class hamil 2}
\cH (\bx_1,\ldots,\bx_N) = \frac{\pi}{2} \sum_{j=1} ^N |\bx_j| ^2 - \sum_{1\leq i < j \leq N} \log |\bx_i - \bx_j| + \cW (\bx_1,\ldots,\bx_N) 
\end{equation}
with $\bx_1,\ldots,\bx_N \in \R ^2$ and $\cW$ a (quite possibly $N$-dependent) function superharmonic in each variable\footnote{Such functions are sometimes called ``plurisuperharmonic'' in the literature. Elliott suggested that we stick to ``superharmonic in each variable'' when writing~\cite{LieRouYng-16,LieRouYng-17}, on the grounds that this was 2016, and that in the Trump era, simple words should be preferred.}:
\begin{equation}\label{eq:superharm 2}
 -\Delta_{\bx_j} \cW \geq 0, \quad \forall j. 
\end{equation}
The first term in~\eqref{eq:class hamil 2} has a constant Laplacian, hence corresponds to the potential generated by a constant neutralizing background, of density $1$ in our units. We proved in~\cite{LieRouYng-16,LieRouYng-17} that the density of charge in the ground state cannot exceed that of the background, on any length scale much larger than the typical inter-particle distance (namely, $1$ in these units).

\begin{theorem}[\textbf{Incompressibility for 2D Coulomb ground states}]\label{thm:incomp GS}\mbox{}\\
There exists a bounded function $g:\R^+ \mapsto \R^+$, independent of $N$ and $\cW$, with 
$$
g(R) \underset{R\to \infty}{\to} 0,
$$ 
such that, for any $X_N^0 = (\bx_1^0,\ldots,\bx_N ^0)$ minimizing $\cH$, any point $a\in \R^2$ and any radius $R > 0$
\begin{equation}\label{eq:incomp GS}
N(a,R) := \sharp\left\{ \bx_j^0 \in X_N ^0 \cap D(a,R) \right\} \leq \pi R^2 (1 + g(R)) 
\end{equation}
where $D(a,R)$ is the disk of center $a$ and radius $R$ and $\sharp$ stands for the cardinal of a discrete set. 
\end{theorem}

A first observation is that, due to the superharmonicity~\eqref{eq:superharm 2}, the proof of Theorem~\ref{thm:sep} applies to this system (there is simply one more superharmonic term in~\eqref{eq:Lieb}), and shows that the minimal distance between points is any case larger than $1/\sqrt{\pi}$. Hence, one can place a disk of radius $1/(2\sqrt{\pi})$ around each point without any overlap between the disks. This leads to a non-trivial bound on the density, but $4$ times too large, something we used in~\cite{RouYng-15} to obtain~\eqref{eq:incomp vague} with $2B\pi \ell^{-1}$ in the right-hand side. 

To obtain the optimal bound needed as an input for Theorem~\ref{thm:ener}, a new idea is needed. In the proof of Theorem~\ref{thm:sep}, we used that any point neutralizes the background in a disk of radius $1/\sqrt{\pi}$ around it (i.e. the total charge of the disk generates no field in its exterior). This can be called a ``screening region'' for a single point charge, and the main tool in the proof of Theorem~\ref{thm:incomp GS} is to define such a screening region associated to \emph{any} discrete set of point charges. This idea was known in potential theory under other names prior to our work, see Remark~\ref{rem:balayage} below.

\subsection{Screening regions}

For neutrality the charge contained in a screening region must be equal to the number of points in the region. Since this is a necessary condition, we include it in the 

\begin{definition}[\textbf{Screening regions}]\label{def:screening}\mbox{}\\
 Let $\bx_1,\ldots,\bx_K$ be points in $\R^2$ and $\Sigma$ an open set with Lebesgue measure 
\begin{equation}\label{eq:size}
 \left| \Sigma\right| = K 
\end{equation}
We say that $\Sigma  = \Sigma(\bx_1,\ldots,\bx_K) \subset \R ^2$ is a \emph{screening region} for the points $\bx_1,\ldots,\bx_K$ if the total electrostatic potential 
\begin{equation}\label{eq:potential}
 \Phi := -\log |\,.\,| \star \left( \sum_{k=1} ^K \delta_{\bx_j} - \one_{\Sigma} \right) 
\end{equation}
satisfies 
\begin{equation}\label{eq:PhiSigma}
\begin{cases}
\Phi &> 0 \mbox{ almost everywhere in } \Sigma\\  
\Phi &= 0 \mbox{ almost everywhere in the complement of } \Sigma. 
\end{cases}
\end{equation}
\end{definition}

There is a saying\footnote{Attributed to James Glimm in one of the citations opening the chapters of the Reed-Simon book series, if my memory is correct.} that ``a good definition is the assumption of a theorem''. At the very least, such definitions are a noteworthy subset of all good definitions. Somewhat dually, a noteworthy subset of all good theorems are ``theorems which prove that a natural definition is not empty''. The following belongs to this class:

\begin{theorem}[\textbf{Incompressible Thomas-Fermi molecules and screening regions}]\label{thm:TF theory}\mbox{}\\
Let $\bx_1,\ldots,\bx_K$ be points in $\R^2$. Consider the energy functional 
\begin{equation}\label{eq:TF func} 
\mathcal \ETF [\sigma]=\int_{\R^2} \sum_{i=1}^K \log|\bx-\bx_i|\sigma(\bx)\,d\bx -\frac{1}{2} \iint_{\R ^2 \times \R ^2} \sigma(\bx)\log|\bx-\by|\sigma(\by)\, d\bx\, d\by. 
\end{equation}
It has a unique minimizer  $\TFmin$  in the class 
\begin{equation}\label{eq:MTF}
\TFM := \left\{\sigma\in L^\infty(\mathbb R^2)\cap  L^1\left(\mathbb R^2, \log(2+|\bx|)d\bx\right), 0\leq \sigma \leq 1, \, \int_{\R ^2}\sigma = K\right\}.
\end{equation}
Moreover $\TFmin$ is of the form 
$$ \TFmin = \1_{\Sigma}$$
for an open set which is a \emph{screening region} for the points $\bx_1,\ldots,\bx_K$. 

In addition, let 
\begin{equation}\label{eq:TF pot}
  \TFpot :=  \TFmin \star \log |\,.\,| - \sum_{i=1}^K \log|\bx-\bx_i|. 
 \end{equation}
For any $R>\max |\bx_i|$ we have
\begin{equation}\label{eq:support TF}
\Sigma \subset D\left(0, R + \sqrt{M_R}\right)
\end{equation}
with
 \beq\label{eq:MR} M_R := \frac{1}{\pi}\sup_{|\bx| = R} \left| \TFpot (\bx) \right|.\eeq
\end{theorem}

The functional minimized to obtain the screening region is referred to as ``incompressible Thomas-Fermi'', for it is reminiscent of a semi-classical approximation of the energy of a 2D molecule (with ``nuclei'' at $\bx_1,\ldots,\bx_K$ and a continuous density of ``electrons'' $\sigma$). We use the word ``incompressible'' because we impose the constraint $\sigma \leq 1$. A ``real'' 2D Thomas-Fermi theory would instead have a penalizing term
$$ \int_{\R^2} \tau(\sigma)$$
in the energy, with $ \tau(\sigma) = 2\pi \sigma^2$, the semi-classical energy density of the free (quantum) 2D electron gas at density $\sigma$. The minimization problem \eqref{eq:TF func}--\eqref{eq:MTF}, corresponds formally to taking $\tau(\sigma)=\sigma^p$ with $p\to\infty$ to enforce the uniform upper bound $\sigma\leq 1$.

The bound~\eqref{eq:support TF} will be important in the sequel, for it gives a control on the shape of the screening region, which could in general be somewhat wild.

\begin{remark}[Subharmonic quadrature domains and partial balayage]\label{rem:balayage}\mbox{}\\
When proof-reading this text, I became aware of the fact (unbeknownst to us when working on~\cite{LieRouYng-16,LieRouYng-17}) that what we call ``screening regions'' in Definition~\ref{def:screening} were known as ``subharmonic quadrature domains'' in the potential theory literature~\cite{GusSha-95,GusPut-07,Sakai-82}. The method to construct such sets in Theorem~\ref{thm:TF theory} is itself known as ``partial balayage of the measure $\sum_{k=1} ^K \delta_{\bx_j}$ to the Lebesgue measure''~\cite{Gustafsson-02}. Another name used more on the physics side~\cite{Zidarov-90} is ``equigravitational mass scattering''. Most of the content of Theorem~\ref{thm:TF theory} can be found in a variety of sources~\cite{Sakai-82,Gustafsson-90,GusSak-94,GarSjo-09,Sjodin-07}, see also~\cite{BalHar-09,GusRoo-18} and references therein. All this is in turn connected to the classical obstacle problem. Our method of proof seems to differ from those in these references, although there is certainly some overlap that we were unaware of.\hfill$\diamond$
\end{remark}

\subsection{Exclusion by screening}

Equipped with the above concept we can now briefly sketch the proof of Theorem~\ref{thm:incomp GS}. Actually we sketch one of two proofs presented in~\cite{LieRouYng-17}, relying on~\eqref{eq:support TF}. The other proof we provided is similar in its first steps, but does not use~\eqref{eq:support TF}. 

\medskip 

\noindent \textbf{Step 1.} Consider a minimizing configuration, a subset thereof, and the screening region associated to it by Theorem~\ref{thm:TF theory}. Since $\cW$ in~\eqref{eq:class hamil 2} is superharmonic in each variable, the arguments of the proof of Theorem~\ref{thm:sep} imply that no other point of the configuration can lie in the screening region just defined. We refer to this as the \emph{exclusion rule}. In the sequel we can forget about minimizers of~\eqref{eq:class hamil 2} and consider all point configurations satisfying this rule that no point can lie within the screening region defined by any subset of other points of the configuration.

\medskip 

\noindent \textbf{Step 2.} Consider the set of all configurations satisfying the exclusion rule. We already know (by the argument sketched below Theorem~\ref{thm:incomp GS}, using screening regions for single points) that the density of any such configuration is bounded above by $4$ (the precise sense in which this holds is similar to~\eqref{eq:incomp GS}, but with a $4$ in the right-hand side). Hence the maximal (or supremal) density $\rho_{\rm max}$ achievable by a configuration  satisfying the exclusion rule is a well-defined number, and we aim at proving that $\rho_{\rm max} \leq 1$.

\medskip 

\noindent \textbf{Step 3.} Consider now a configuration $(\bx_1,\ldots,\bx_k,\ldots)$ satisfying the exclusion rule, and achieving the maximal density $\rho_{\rm max}$. This configuration cannot have any large vacancy. Indeed, a density lower than $\rho_{\rm max}$ in some region would have to be compensated for by a density higher than $\rho_{\rm max}$ in another region. This would contradict the definition of $\rho_{\rm max}$ as the maximal density achievable while satisfying the exclusion rule. 

\medskip 

\begin{figure}[h]
\begin{center}
\includegraphics[width=10cm]{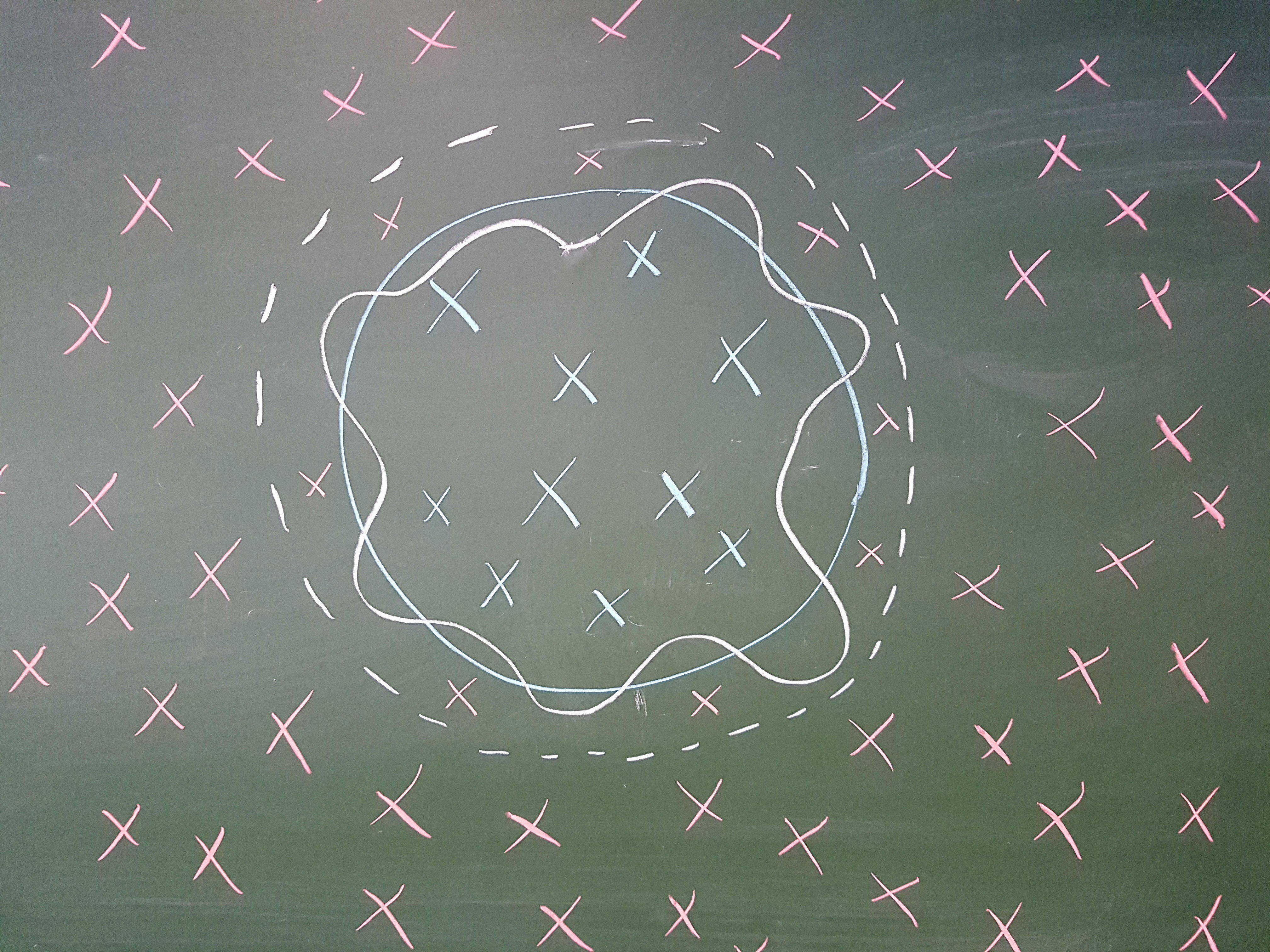}\\
%
%
%
\caption{\textbf{Good configuration}. The blue points inside the (blue) circle generate a screening region (inside of white line), avoiding the other (red) points. It is contained in a disk (inside of white dashed circle) not too large compared to the original blue circle.}
\label{fig:good}
\end{center}
 \end{figure}

\noindent \textbf{Main and final step.} It is hence sufficient (roughly) to consider a configuration satisfying the exclusion rule and having a maximal density $\rho_{\rm max}$ ``everywhere''. Consider the points in said configuration lying in some disk of radius $R$ (say, without loss, of center $0$), and the associated screening region $\Sigma_R$. Since, by definition 
$$ |\Sigma_R| = \mbox{ number of points in }  D(0,R)$$
the conclusion of the theorem will follow if we prove that 
$$ |\Sigma_R| \lesssim \pi R^2$$
for large $R$. This we do by proving that $\Sigma_R$ cannot ``leak'' too much out of the original disk $D(0,R)$. Namely, we want to be in the ``good case'' sketched in Figure~\ref{fig:good}, where the screening region is included in a slightly larger disk of radius $R+r$, $r\ll R$.

\begin{figure}[h]
\begin{center}
\includegraphics[width=10cm]{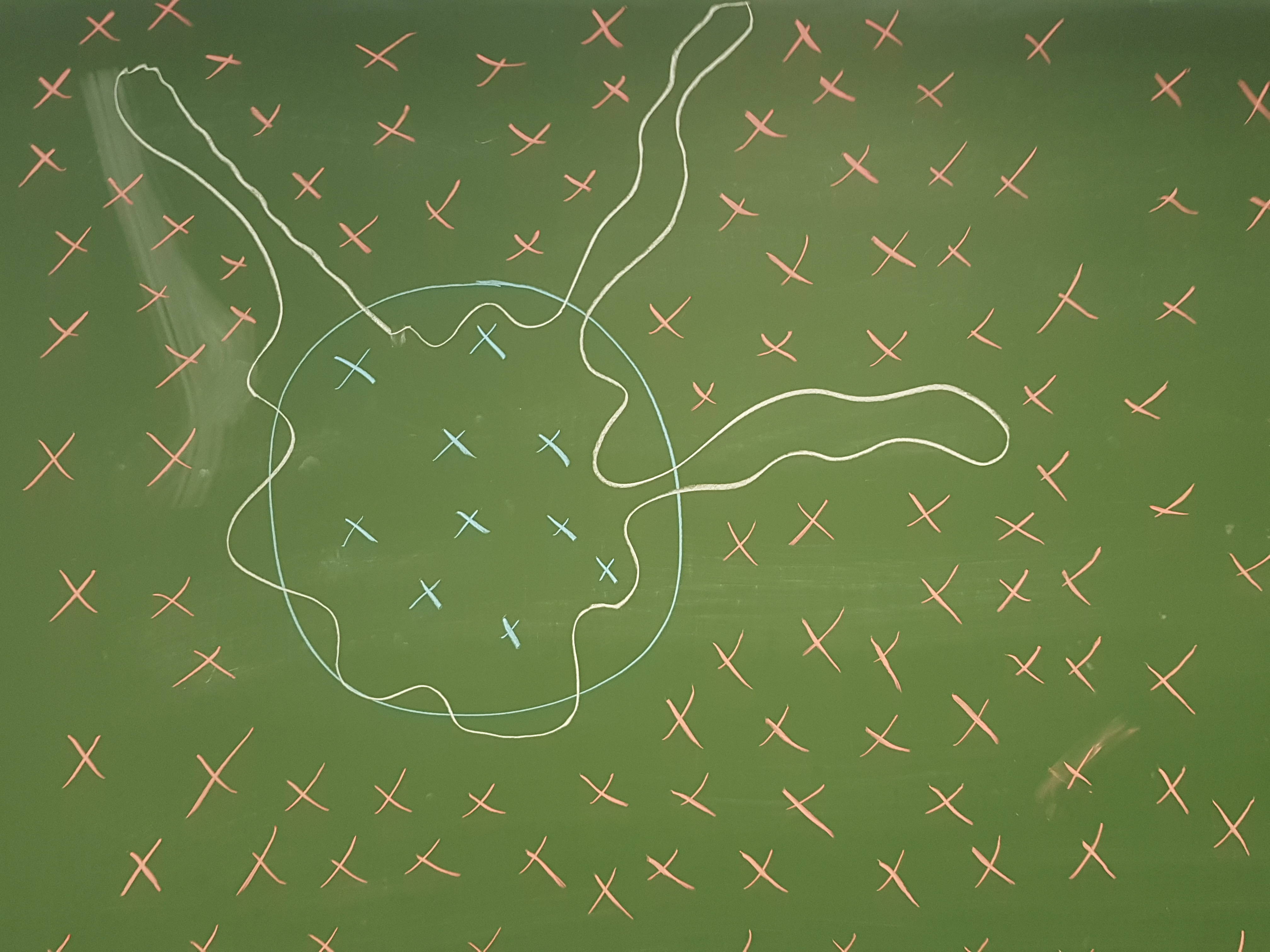}\\
%
%
%
\caption{\textbf{Pathological configuration, to be excluded}. The blue points inside the (blue) circle generate a screening region (inside of white line). It sends tendrils out to infinity, while still avoiding the other (red) points.}
\label{fig:bad}
\end{center}
 \end{figure}

Since the configuration satisfies the exclusion rule, the screening region has to avoid all the points in the exterior of the disk. Our main enemy is thus, as sketched in Figure~\ref{fig:bad} that a tendril of the screening region is sent to infinity, winding its way bizarrely around all the other points.

Let $\Phi$ be the potential~\eqref{eq:potential} generated by the points in the original disk and the screening region. By the exclusion rule and~\eqref{eq:PhiSigma}, it must vanish at all the configuration's points outside of the disk. Since the configuration may not have large vacancies, the points the screening region must avoid, at which $\Phi = 0$, are numerous. A few estimates prove that these points are sufficiently dense to deduce 
$$  \sup_{\bx \in \partial D(0,R)} \left| \Phi (\bx) \right| \ll R^2.$$
Using~\eqref{eq:support TF} we deduce that 
$$ \Sigma_R \subset D(0,R+r)$$
with $r\ll R$, and hence that 
\begin{equation}\label{eq:concl}
  \mbox{ number of points in }  D(0,R) = |\Sigma_R| \leq \pi R^2 (1+o_{R\to \infty}(1)), 
\end{equation}
the desired conclusion.

\begin{remark}[A proof variant]\label{rem:variant}\mbox{}\\
According to~\cite[Theorem~5.4]{Gustafsson-02}, $\Sigma_R$ can be written as a union of disks with centers in $D(0,R)$ (the proof is in~\cite{GusSak-03,GusSak-04}). This also excludes the pathological configuration in Figure~\ref{fig:bad}. Indeed, write 
$$ \Sigma_R = \bigcup_{a \in \Sigma \cap D(0,R)} D(a,r(a))$$
with $r(a)>0$. Let $\by$ be the point in $\Sigma_R$ the furthest away from $D(0,R)$. Then $\by \in D(a,r(a))$ for some $a\in D(0,R)\cap \Sigma_R$ and $D(a,r(a)) \subset \Sigma_R$, so that $R + r = R + r(a) \geq |\by|$. The part of the disk $D(a,r)$ which is not included in $D(0,R)$ has size $\sim r^2$, so if $r\to \infty$ when $R\to \infty$, it must touch points from the configuration outside $D(0,R)$ (recall that the configuration may not have large vacancies). This would be a contradiction, so $r$ must stay bounded by a (possibly large) constant when $R\to \infty$. It follows that $|\by| \leq C$, i.e.
$$ \Sigma_R \subset D (0,R+C)$$
for some constant $C$, and thus we get~\eqref{eq:concl} again.
 \hfill$\diamond$
\end{remark}

\section{Short Conclusion}

We have discussed an unusual variational problem in many-body quantum mechanics, and the motivation for introducing it, which comes from fractional quantum Hall (FQH) physics. The aim is to show that the Laughlin state with quasi-holes is stable under weak external and interaction potentials. The approach to this problem we have been following in the past few years (with Elliott Lieb, Alessandro Olgiati, Sylvia Serfaty and Jakob Yngvason) proceeds by analogy with the study of (somewhat contrived, if seen independently from the FQH motivation) classical Gibbs states of Coulomb systems. This lead us to a brief and partial review of known results bearing on more standard classical Coulomb systems: the homogeneous and inhomogeneous jellium. An unpublished (but generously communicated to people that had the use for it) theorem of Elliott Lieb bearing on such systems then provided the source of inspiration for the derivation of the main tools used in the study of the FQH variational problem.

\medskip

\textbf{Acknowledgments:} I am financially supported by the European Research Council (under the European Union's Horizon 2020 Research and Innovation Programme, Grant agreement CORFRONMAT No 758620). It is a pleasure to thank the editors of the present collection for inviting me to write this text, and the aforementioned collaborators, work with whom it is based on. I am indebted to Bj\"orn Gustafsson for useful discussions and references related to Remarks~\ref{rem:balayage} and~\ref{rem:variant}. Finally it is an honor to thank Elliott H. Lieb for the opening or widening of so many fields of mathematical physics, for us to continue exploring. 

\newpage
%

\end{document}